\documentclass[10pt,a4paper,reqno]{article}
\usepackage{geometry}

\usepackage[pdftitle={Holomorphic Poisson Manifolds and Holomorphic Lie Algebroids},pdfauthor={Camille Laurent, Mathieu Stienon, 
Ping Xu}]{hyperref}

\usepackage[T1]{fontenc}
\usepackage[latin1]{inputenc}
\usepackage{lmodern}

\usepackage{setspace}
\usepackage{indentfirst}

\usepackage[intlimits]{amsmath}
\usepackage{amssymb}
\usepackage{stmaryrd}

\usepackage{graphicx}
\usepackage[all]{xy}
\def\dar[#1]{\ar@<2pt>[#1]\ar@<-2pt>[#1]}
\usepackage{amsthm}

\theoremstyle{plain}% default
\newtheorem{prop}{Proposition}[section]
\newtheorem{lem}[prop]{Lemma}

\newtheorem{cor}[prop]{Corollary}
\newtheorem{thm}[prop]{Theorem}

\newtheorem*{prop*}{Proposition}
\newtheorem*{lem*}{Lemma}
\newtheorem*{sublem*}{Sublemma}
\newtheorem*{cor*}{Corollary}
\newtheorem*{thm*}{Theorem}
\newtheorem*{hypo*}{Hypothesis}
\newtheorem*{question*}{Question}
\newtheorem*{conjecture*}{Conjecture}
\newtheorem*{scholum*}{Scholum}
\newtheorem{defn}[prop]{Definition}
\newtheorem*{defn*}{Definition}

\theoremstyle{definition}

\newtheorem*{con*}{Construction}
\newtheorem*{note*}{Note}

\theoremstyle{remark}

\newtheorem*{warning*}{Warning}
\newtheorem*{shortnote*}{Note}
\newtheorem*{claim*}{Claim}
\newtheorem*{axiom*}{Axiom}

\newtheoremstyle{slanted}% name
  {3pt}%      Space above, empty = `usual value'
  {3pt}%      Space below
  {\slshape}% Body font
  {}%         Indent amount (empty = no indent, \parindent = para indent)
  {\bfseries}% Thm head font
  {.}%        Punctuation after thm head
  {.5em}%     Space after thm head: " " = normal interword space;
  {}% Thm head spec

\theoremstyle{slanted}

\newtheorem*{example*}{Example}
\newtheorem*{examples*}{Examples}
\newtheorem{ex}[prop]{Example}

\newtheorem*{ex*}{Example}
\newtheorem*{exs*}{Examples}

\newtheorem*{remark*}{Remark}
\newtheorem*{remarks*}{Remarks}
\newtheorem{rmk}[prop]{Remark}

\newtheorem*{rmk*}{Remark}
\newtheorem*{rmks*}{Remarks}

\DeclareMathOperator{\idn}{id}

\DeclareMathOperator{\Hom}{Hom}

\DeclareMathOperator{\pr}{pr}

\DeclareMathOperator{\DER}{Der}

\newcommand{\beq}[1]{\begin{equation}\label{#1}}
\newcommand{\eeq}{\end{equation}}

\newcommand{\JJ}{\mathbb{J}} % generalized complex structure in this paper

\newcommand{\CC}{\mathbb{C}}
\newcommand{\RR}{\mathbb{R}}

\newcommand{\derlie}[1]{\mathcal{L}_{#1}}

\newcommand{\ii}{\mathbin{\vrule width1.5ex height.4pt\vrule height1.5ex}}
\newcommand{\per}{\ii}
\newcommand{\pairing}[2]{\langle #1,#2 \rangle} % inner product
\newcommand{\lie}[2]{[#1,#2]} % Lie bracket
\newcommand{\schouten}[2]{[#1,#2]} % Schouten bracket
\newcommand{\courant}[2]{\llbracket#1,#2\rrbracket} % Courant bracket
\newcommand{\anchor}{\rho} % anchor map
\newcommand{\delbar}{\bar{\partial}}
\newcommand{\del}{\partial}

\newcommand{\thalf}{\tfrac{1}{2}}

\newcommand{\rond}{\circ}

\newcommand{\cc}[1]{\overline{#1}} % complex conjugate of

\newcommand{\gendex}[2]{\left\{ #1 \right\}_{#2}}
\newcommand{\genrel}[2]{\left\{ #1 | #2 \right\}}

\newcommand{\cinf}[1]{C^{\infty}(#1)}
\newcommand{\sections}[1]{\Gamma(#1)}
\newcommand{\XX}{\mathfrak{X}} % vector fields
\newcommand{\vf}{\mathfrak{X}} % vector fields
\newcommand{\df}{\Omega} % differential forms
\newcommand{\injection}{\hookrightarrow}

\newcommand{\isomorphism}{\cong}%{\simeq}

\newcommand{\ssi}{\Leftrightarrow}

\newcommand{\diese}{^{\sharp}}

\newcommand{\cmplx}{_{\mathbb{C}}} % complexified
\newcommand{\inv}{^{-1}}

\newcommand{\TTM}{TM\oplus T^*M}

\newcommand{\dpi}{d_\pi}

\newcommand{\NIJ}{\phi} % de Rham
\newcommand{\diff}{{\rm d}} % de Rham
\newcommand{\gm}{\Gamma}

\newcommand{\hf}[1]{\mathcal{O}_{#1}} % sheaf of holomorphic functions on
\newcommand{\shs}[1]{\mathcal{#1}} % stands for sheaf of holomorphic sections of holomorphic bundle
\newcommand{\pire}{\pi_R} % real part of bivector \pi
\newcommand{\piim}{\pi_I} % imaginary part of bivector \pi
\newcommand{\picre}{\pi_\Re} % real part of bivector \pi
\newcommand{\picim}{\pi_\Im} % imaginary part of bivector \pi
\newcommand{\pb}[2]{\{#1,#2\}} % Poisson bracket on functions associated to holomorphic \pi
\newcommand{\pbre}[2]{\{#1,#2\}_R} % Poisson bracket on functions associated to \pire
\newcommand{\pbim}[2]{\{#1,#2\}_I} % Poisson bracket on functions associated to \piim
\newcommand{\pbcre}[2]{\{#1,#2\}_\Re} % Poisson bracket on functions associated to \pire
\newcommand{\pbcim}[2]{\{#1,#2\}_\Im} % Poisson bracket on functions associated to \piim
\newcommand{\com}{\bullet} % Kai's bullet
\newcommand{\RA}[1]{\Hom_\RR (#1,\RR)}% real dual of  A
\newcommand{\CA}[1]{ \Hom_\CC (#1,\CC)}% complex dual of A
\newcommand{\be }{\begin{eqnarray*}}
\newcommand{\ee }{\end{eqnarray*}}
\newcommand{\canchor}{\anchor_\CC}

\newcommand{\xii}{\xi^{1,0}}
\newcommand{\etaa}{\eta^{1,0}}

\newcommand{\OA}[2]{\Omega^{0, #1}_X\otimes_{C^\infty_X}\shs{A}^{#2, 0}_\infty}
\newcommand{\ppi}{\hat{\pi}}

\newcommand{\fboxr}[1]{#1}

\newcommand{\arxiv}[1]{\href{http://arxiv.org/abs/#1}{\texttt{arXiv:#1}}}

\begin{document}

\title{Holomorphic Poisson  Manifolds and Holomorphic Lie Algebroids}
\author{
Camille Laurent-Gengoux \\ 
D\'epartement de math\'ematiques \\ Universit\'e de Poitiers \\ 
86962 Futuroscope-Chasseneuil, France \\ 
\href{mailto:laurent@math.univ-poitiers.fr}{\texttt{laurent@math.univ-poitiers.fr}}
\and 
Mathieu Sti\'enon
\thanks{Research supported by the European Union through the FP6 Marie Curie R.T.N. ENIGMA
(Contract number MRTN-CT-2004-5652).} \\
E.T.H.~Z\"urich \\ Departement Mathematik \\ 
8092 Z\"urich, Switzerland \\
\href{mailto:stienon@math.ethz.ch}{\texttt{stienon@math.ethz.ch}}
\and 
Ping Xu
\thanks{Research partially supported by NSF
grants DMS-0306665 and DMS-0605725 \&  NSA grant H98230-06-1-0047} \\
Department of Mathematics \\ Penn State University \\ 
University Park, PA 16802, U.S.A. \\
\href{mailto:ping@math.psu.edu}{\texttt{ping@math.psu.edu}} }
%\date{\texttt{\jobname.tex}}
\date{}
%\subjclass{}
\maketitle

\centerline{{Dedicated to the memory of Paulette Libermann}}

\begin{abstract}
We study holomorphic Poisson manifolds and  holomorphic Lie algebroids 
 from the viewpoint of real Poisson
geometry. We give a characterization of holomorphic Poisson 
structures in terms of the Poisson Nijenhuis structures of
 Magri-Morosi and describe a double complex which computes the
 holomorphic Poisson cohomology. A holomorphic Lie algebroid structure
on a vector bundle $A\to X$ is shown to be equivalent to
 a matched pair of complex Lie algebroids $(T^{0,1}X,A^{1,0})$, 
 in the sense of Lu.
The holomorphic Lie algebroid cohomology of $A$ is isomorphic
to the cohomology of the elliptic Lie algebroid
$T^{0,1}X\bowtie A^{1,0}$. 
In  the case when $(X,\pi)$ is a holomorphic
Poisson manifold and $A=(T^*X)_\pi$,
such an elliptic  Lie algebroid coincides with the Dirac
structure corresponding to the associated
generalized complex structure of the holomorphic
Poisson manifold.
\end{abstract}

\newpage
\tableofcontents

\section{Introduction}

The aim of this paper is to solve several problems naturally arisen in studying the connection between holomorphic Poisson manifolds and holomorphic Lie algebroids with  real Poisson geometry.

Holomorphic Poisson structures appear naturally in many places \cite{MR1665693,MR1396600,MR1178029,MR1449222,
MR1980616,MR2104604,MR2263715,MR2372206,MR2245536}.
For instance, any semi-simple complex Lie group admits
a natural Poisson group structure \cite{MR688240, MR725930,MR842417},
 which is holomorphic. 
Its dual is also a holomorphic Poisson group.
Indeed one of the simplest types of examples of holomorphic Poisson manifolds
are  the Lie-Poisson structures on the dual of  complex Lie algebras.
Holomorphic Poisson structures were also studied from
the point of view of algebraic geometry by Bondal \cite{Bondal1993}
and Polishchuk \cite{MR1465521} in the middle of the 90's.
Recently, holomorphic Poisson structures were linked to 
generalized complex geometry \cite{Gualtieri2007a, Gualtieri2007,MR2013140,MR2217300,MR2055289}.

There are several equivalent ways of defining holomorphic Poisson structures.
One simple definition is, like in the real case, a holomorphic
 bivector field $\pi$ (i.e. $\pi\in \gm (\wedge^2 T^{1, 0}X)$
such that $\delbar\pi=0$) satisfying the equation $\schouten{\pi}{\pi}=0$.
Since $\wedge^2 T_{\CC}X = \wedge^2 TX \oplus i \wedge^2 TX$,
for any $\pi \in \sections{\wedge^2 T_{\CC}X}$, we can write
$\pi=\pire+i\piim$, where $\pire$ and $\piim\in\sections{\wedge^2 TX}$
are  bivector fields on the underlying real manifold $X$.

\begin{quote}
\textbf{Problem 1.}
Are $\pire$ and $\piim$ Poisson structures? And conversely, given two
Poisson structures $\pire$ and $\piim$, when does $\pi=\pire+i\piim$
define a holomorphic Poisson structure?
\end{quote} 

We give an affirmative answer to the first question. As for the second, we show that
$\pi=\pire+i\piim$ is holomorphic Poisson if, and only if, 
$(\piim,J)$ is a Poisson Nijenhuis structure and
$\pire\diese=J\rond\piim\diese$. Thus   $(\pire,\piim)$ 
is a bi-Hamiltonian structure on $X$.

Poisson Nijenhuis structures were introduced
by Magri and Morosi \cite{MR773513, MR900387} in their study of
bi-Hamiltonian systems, and were intensively studied afterwards \cite{MR1077465,MR1390832}.
A Poisson Nijenhuis structure \cite{MR1077465,MR1421686} on a 
manifold $X$ consists of a pair $(\pi, N)$, where $\pi$ is a Poisson tensor
on $X$ and $N : TX\to TX$ is a Nijenhuis tensor,
which satisfy some  compatibility conditions (see Section~\ref{section:2.3} for the precise conditions).
By a Nijenhuis tensor, we mean a 
$(1,1)$-tensor on $X$ with vanishing Nijenhuis torsion.

Since Poisson Nijenhuis structures
are related to generalized complex structures \cite{Crainic2007,MR2276462},
as a consequence, we recover the well known correspondence between holomorphic Poisson structures and generalized complex structures (of a special type) \cite{MR2217300,Gualtieri2007a}.

Another natural question is:

\begin{quote}
\textbf{Problem 2.}
Given a holomorphic Poisson structure $\pi=\pire+i\piim$,
are the holomorphic symplectic foliation of $\pi$ and the symplectic
foliations of $\pire$ and $\piim$ related?
\end{quote}

Indeed we show that all these symplectic foliations 
coincide. Also for a holomorphic symplectic
$2$-form $\omega=\omega_R+i\omega_I$, we show that
the real and imaginary parts of its
holomorphic Poisson tensor are, up to
a constant scalar, the Poisson tensors corresponding
to $\omega_R$ and $\omega_I$, respectively.

Lie algebroids are an extremely powerful tool in Poisson geometry. Indeed the Lie algebroid structures on a given vector bundle are in one-one correspondence with the so called fiberwise linear Poisson structures on the dual bundle. 
This correspondence extends to the holomorphic context; 
any holomorphic Lie algebroid structure on the vector bundle $A\to X$ gives rise to a fiberwise linear holomorphic Poisson structure on $A^*\to X$. Thus the real and imaginary
parts of this holomorphic Poisson structure
are fiberwise linear Poisson structures on the dual bundle (being
considered as a real vector bundle). Hence one obtains
two real Lie algebroid structures  $A_\Re$ and $A_\Im$, respectively.

\begin{quote}
\textbf{Problem 3.}
Obtain an explicit description of the Lie algebroid structures  $A_\Re$ and $A_\Im$ in terms of the holomorphic Lie algebroid structure on $A$.
\end{quote}

Let $A\to X$ be a vector bundle endowed with a holomorphic Lie algebroid structure. Extending the Lie bracket on the space of holomorphic sections of $A\to X$ to the space of all smooth sections so as to preserve the Leibniz rule, we get a real Lie algebroid structure $A_R$ on the bundle $A\to X$. It turns out that, up to a scalar constant, $A_\Re$ is isomorphic to $A_R$. 
The multiplication by $\sqrt{-1}$ in the fibers of 
$A$ defines a real vector bundle map $j: A_R\to A_R$ over the
identity map satisfying $j^2=-\idn$. 

We prove that the Nijenhuis torsion of $j:A_R\to A_R$ 
vanishes and that, up to a scalar constant, $A_\Im$ is 
isomorphic to $(A_R)_j$, the deformation of the Lie algebroid $A_R$ by $j$. 
Extending $j$ by  $\CC$-linearity, we get a bundle map $j: A_\CC\to A_\CC$ with $j^2=-\idn$.
Since the Nijenhuis torsion of $j$ vanishes, its eigenbundles $A^{1,0}$ and $A^{0,1}$ with eigenvalues $i$ and $-i$ are complex Lie algebroids.

There is yet another connection between Poisson manifolds
and Lie algebroids. Given a Poisson manifold $(X,\pi)$,
it is well known that $T^*X$ carries a natural Lie algebroid structure.  This holds for holomorphic Poisson structures as well. Namely, if 
$(X,\pi)$ is  a holomorphic Poisson manifold, then $T^*X$ is naturally a 
holomorphic Lie algebroid, denoted $(T^*X)_\pi$.
On the other hand, as highlighted earlier, each
holomorphic Poisson structure corresponds to a 
generalized complex structure $\JJ$, whose $(-i)$-eigenbundle $L$ is a Dirac structure and thus a complex Lie algebroid \cite{MR2285039,MR1718638}.

\begin{quote}
\textbf{Problem 4.}
What is the precise relation between the holomorphic Lie algebroid
$(T^*X)_\pi$ and the complex Lie algebroid $L$?
\end{quote}

A key ingredient to answer this problem is the notion of
matched pairs studied by Lu, Mackenzie, and Mokri \cite{MR1430434,Mackenzie2007,MR1460632}.
We show that $(T^{0,1}X,(T^{1,0}X)^*_\pi)$ is a matched
pair (here $(T^{1, 0}X)^*_\pi=A^{1,0}$ for $A=(T^*X)_\pi$)
 and  $T^{0,1}X\bowtie(T^{1,0}X)^*_\pi$
(see Theorem~\ref{thm:5.2} for the definition of $\bowtie$) 
is isomorphic to $L$. Furthermore, we prove that the holomorphic Poisson cohomology of $\pi$, which is defined to be the holomorphic Lie algebroid cohomology of $(T^*X)_\pi$, is isomorphic to the cohomology 
of the elliptic Lie algebroid $T^{0,1}X\bowtie(T^{1,0}X)^*_\pi$.
This leads to our next problem:

\begin{quote}
\textbf{Problem 5.}
Given an arbitrary holomorphic Lie algebroid $A$, find a complex Lie algebroid $L$ whose cohomology groups are isomorphic to those of $A$.
\end{quote}

The cohomology of a holomorphic Lie algebroid $A$ 
is the cohomology of the complex of sheaves 
$(\Omega_A^\bullet,d_A)$ as introduced by Evens-Lu-Weinstein \cite{MR1726784}
 (see Definition~\ref{defn:5.11}),
while the cohomology of a complex (smooth) 
Lie algebroid $L$ is the cohomology of the cochain complex 
$(\sections{\wedge^\bullet L^*},d_L)$. 
So, in a certain sense, solving the problem above amounts to
finding a Dolbeault type of resolution for arbitrary holomorphic
Lie algebroids. 

The solution is $L=T^{0,1}X\bowtie A^{1,0}$.
Indeed we show that $A$ is holomorphic if, and only if, 
$(T^{0,1}X,A^{1,0})$ is a matched pair (see Theorem~\ref{HLA=MP});
one may thus form the complex Lie algebroid 
$T^{0,1}X\bowtie A^{1,0}$, which is in fact an elliptic Lie algebroid
in the sense of Block \cite{Block2007}. 
The Lie algebroid cohomology of the latter can be expressed 
as the total cohomology of a double complex.

The following notations are widely used in the sequel. 
For a manifold $M$, we use $q_M$ to denote the projection $TM\to M$. 
And given a complex manifold $X$, $T_\CC X$ is shorthand for the complexified tangent bundle $TX\otimes\CC$ while $T^{1,0} X$ (resp. $T^{0,1} X$) stands for the $+i$- (resp. $-i$-) eigenbundle of the almost complex structure. By $\XX^{k,l}(X)$ we denote the
space of sections of $\wedge^k T^{1,0}X\otimes \wedge^l T^{0,1}X\to X$,
and by $\Omega^{k,l}(X)$ the space of differential forms of type $(k,l)$.
For a Lie algebroid $A$, 
the Nijenhuis torsion \cite{MR1077465,MR1421686} of a bundle map
$\NIJ: A\to A$ over the identity is  denoted 
$\mathcal{N}_\NIJ$, which is a section in
$\sections{\wedge^2 A^*\otimes A}$ defined by
\begin{equation}
\label{eq:nij_def}
\mathcal{N}_\NIJ(V,W)
=\lie{\NIJ V}{\NIJ W}-\NIJ(\lie{\NIJ V}{W}+\lie{V}{\NIJ W}-\NIJ \lie{V}{W}),
\qquad \forall V, W \in \sections{A} .
\end{equation}
When $A$ is the Lie algebroid $TX$ and $\NIJ: TX\to TX$ is
a $(1,1)$-tensor, the Nijenhuis torsion $\mathcal{N}_\NIJ$
is a $(2,1)$-tensor on $X$.

Note that the modular classes of holomorphic Lie algebroids
were studied by Evens-Lu-Weinstein \cite{MR1726784} and
Huebschmann \cite{MR1696093} while the modular classes of holomorphic
Poisson manifolds were studied by Brylinski-Zuckerman \cite{MR1665693}.
In a separate paper, we will investigate the relation between 
these modular classes and their counterparts in real Poisson 
geometry, and in particular with the modular classes of 
Poisson Nijenhuis manifolds recently studied by Damianou-Fernandes 
\cite{Damianou2007} and Kosmann-Schwarzbach-Magri \cite{Kosmann-Schwarzbach2007}.

\paragraph{Acknowledgments}
We would like to thank Centre \'Emile Borel and Peking
University for their hospitality while work on this project was
 being done. We also wish to thank Sam Evens and  Alan Weinstein
for useful discussions and comments.

\section{Holomorphic Poisson manifolds}

\subsection{Definition}

\begin{defn}
A holomorphic Poisson manifold is a complex manifold
$X$ whose sheaf
of holomorphic functions $\hf{X}$ is a sheaf of Poisson algebras.
\end{defn}

By a sheaf of Poisson algebras over $X$,
we mean that, for each open subset
$ U \subset X$, the ring $\hf{X}(U) $ is endowed with a Poisson bracket
such that all restriction maps $\hf{X}(U) \to \hf{X}(V) $ (for arbitrary
open subsets $V \subset U \subset X $) are morphisms of Poisson algebras.
Moreover, given an open subset $U\subset X$, an open covering
$\gendex{U_i}{i\in I}$ of $U$ and a pair of functions $f,g\in\hf{X}(U)$,
then the local data $\pb{f|_{U_i}}{g|_{U_i}}$ $(i\in I)$ glue up to
$\pb{f|_U}{g|_U}$ if they coincide on the overlaps $U_i\cap U_j$.

\begin{lem}
On a given complex manifold $X$, holomorphic
Poisson structures are in one-to-one correspondence with
holomorphic bivector fields $\pi$ (i.e. $\pi\in \gm (\wedge^2 T^{1, 0}X)$
such that $\delbar\pi=0 $), satisfying the equation $\schouten{\pi}{\pi}=0$.
\end{lem}

\begin{proof}
This is a standard  result. For completeness, let us
sketch a proof here. Choose any complex coordinate chart $(U,\phi)$,
which identifies $U\subset X$ with an open ball of  $\CC^{n}$.
As in the smooth case, that  $\hf{X}(U) $  is a
Poisson algebra is equivalent to  the existence of a
holomorphic bivector field
$\pi_U $ on $U$ satisfying the relation $[\pi_U,\pi_U] =0$.
 Moreover the compatibility condition on the restriction maps implies
that there indeed exists a holomorphic bivector field $\pi$ on $X$,
 whose restriction to $U $ is $\pi_U $
for all such open subset $U$.
\end{proof}

\subsection{Associated real Poisson structures}

Since $\wedge^2 T_{\CC}X = \wedge^2 TX \oplus i \wedge^2 TX$,
for any $\pi \in \sections{\wedge^2 T_{\CC}X}$, we can write
$\pi=\pire+i\piim$, where $\pire$ and $\piim\in\sections{\wedge^2 TX}$
are (real) bivector fields on $X $ (seen as a real manifold by
forgetting the complex structure).
Note that sections of $\wedge^2 T_{\mathbb C} X $ (in particular of
$\wedge^2 T^{1,0} X $) can be seen as bidifferential
operators on $C^{\infty}(M,{\mathbb C})$.
The real bivector fields $\pire$ and $\piim$ are then the real
and imaginary parts of these bidifferential operators.

Both $\pire$ and $\piim$ define
brackets $\pbre{\cdot}{\cdot}$ and $\pbim{\cdot}{\cdot}$
on $\cinf{M,\RR}$ in the standard way.
These extend to $\cinf{M,\CC}$ by $\CC$-linearity.
The next lemma describes such an extension.

\begin{lem}\label{lem:1.3}
\begin{enumerate}
\item\label{oqp1}
Under the direct sum decomposition
\[ \wedge^2T\cmplx X=\wedge^2T^{1,0} X\oplus (T^{1,0} X\wedge T^{0,1} X)
\oplus \wedge^2T^{0,1} X, \]
we have
\[ \pire=\tfrac{\pi}{2}+0+\tfrac{\cc{\pi}}{2} \qquad \text{and} \qquad  \piim=\tfrac{\pi}{2i}+0+\tfrac{-\cc{\pi}}{2i}. \]
\item\label{oqp2}
 $\forall f,g\in\hf{X} (U)$, we have the following relations:
\begin{align*}
\pbre{f}{g}&=\tfrac{1}{2}\pb{f}{g}, & \pbre{\cc{f}}{\cc{g}}&=\tfrac{1}{2}\cc{\pb{f}{g}}, & \pbre{f}{\cc{g}}&=0, \\
\pbim{f}{g}&=-\tfrac{i}{2}\pb{f}{g}, & \pbim{\cc{f}}{\cc{g}}&=\tfrac{i}{2}\cc{\pb{f}{g}}, & \pbim{f}{\cc{g}}&=0.
\end{align*}
\item\label{oqp3}  Both $\pire$ and $\piim$ are Poisson tensors.
\end{enumerate}
\end{lem}

\begin{proof} \fboxr{\ref{oqp1}} Immediate.
\fboxr{\ref{oqp2}} If $f,g\in\hf{X} (U)$, then $\cc{f}$ and $\cc{g}$ are antiholomorphic.
Hence $\del\cc{f}=0=\del\cc{g}$. Therefore $\pi(\del f,\del g)=\pb{f}{g}$,
$\pi(\del\cc{f},\del g)=0$ and $\pi(\del\cc{f},\del\cc{g})=0$. The conclusion follows.
\fboxr{\ref{oqp3}} It suffices to prove the Jacobi identity
for $\pbre{\cdot}{\cdot}$ and $\pbim{\cdot}{\cdot}$. 
From \eqref{oqp2}, it follows that it holds for holomorphic functions
in $\hf{X}(U)$. It follows from the Leibniz rule that
the Jacobi identity holds for all complex valued smooth
functions in $\cinf{X,\CC}$. This concludes the proof.
\end{proof}

As an immediate consequence, we have the following

\begin{cor}\label{cor:2.4}
For all $f,g\in\hf{X}(U)$, we have
\begin{align*}
\pbre{\Re f}{\Re g}&=\tfrac{1}{4}\Re\pb{f}{g}, & \pbim{\Re f}{\Re g}&=\tfrac{1}{4}\Im\pb{f}{g}, \\
\pbre{\Im f}{\Im g}&=-\tfrac{1}{4}\Re\pb{f}{g}, & \pbim{\Im f}{\Im g}&=-\tfrac{1}{4}\Im\pb{f}{g}, \\
\pbre{\Re f}{\Im g}&=\tfrac{1}{4}\Im\pb{f}{g}, & \pbim{\Re f}{\Im g}&=-\tfrac{1}{4}\Re\pb{f}{g}.
\end{align*}
where $\Re f$ and $\Im f$ stand for the real and imaginary
part of the function $f$, respectively.

Thus, in a local chart $(z_1=x_1+i y_1,\cdots,z_n=x_n+i y_n)$ of complex coordinates of $X$,
we have
\begin{align*}
\pbre{x_i}{x_j}&=\tfrac{1}{4}\Re\pb{z_i}{z_j}, & \pbim{x_i}{x_j}&=\tfrac{1}{4}\Im\pb{z_i}{z_j}, \\
\pbre{y_i}{y_j}&=-\tfrac{1}{4}\Re\pb{z_i}{z_j}, & \pbim{y_i}{y_j}&=-\tfrac{1}{4}\Im\pb{z_i}{z_j}, \\
\pbre{x_i}{y_j}&=\tfrac{1}{4}\Im\pb{z_i}{z_j}, & \pbim{x_i}{y_j}&=-\tfrac{1}{4}\Re\pb{z_i}{z_j}.
\end{align*}
\end{cor}

\subsection{Poisson Nijenhuis structures}
\label{section:2.3}
\begin{lem}
\label{lem:2.5}
Let $\pi=\pire+i\piim\in\sections{\wedge^2T\cmplx X}$ with $\pire,\piim\in\sections{\wedge^2 TX}$.
Then $\pi\in\sections{\wedge^2 T^{1,0} X}$ iff $\pire\diese=\piim\diese\rond J^*$, where $J: TX\to TX$ is
the almost complex structure on $X$.
\end{lem}

\begin{proof}
We have
\begin{align*}
& \pi\in\sections{\wedge^2 T^{1,0} X} \\
\ssi \qquad &
\pi(\alpha,\beta)=0,\;\forall \alpha \in \Omega^{0, 1}(X), \
\beta\in \Omega^1_\CC (X) \\
\ssi \qquad & \pi(\alpha,\beta)=0,\; \alpha =\tfrac{1+iJ^*}{2} (\alpha '), \ \
\forall \alpha' , \beta\in \Omega^1_\CC (X) \\
\ssi \qquad & \pi\diese\rond\big(\tfrac{1+iJ^*}{2}\big)=0 \\
\ssi \qquad & i \pi\diese=\pi\diese\rond J^* \\
\ssi \qquad & 2\pire\diese\rond J^*=(\pi\diese+\cc{\pi}\diese) \rond
J^*=i(\pi\diese-\cc{\pi}\diese)=-2\piim\diese \\
\ssi \qquad & \pire\diese=\piim\diese\rond J^*
\qedhere \end{align*}
\end{proof}

Recall that a Poisson Nijenhuis structure \cite{MR1077465,MR1421686} on a
manifold $X$  consists of a pair $(\pi, N)$, where $\pi$ is a Poisson tensor
on $X$ and $N : TX\to TX$ is a Nijenhuis tensor
such that the following compatibility
conditions are satisfied:
\begin{gather*}
N\rond\pi\diese=\pi\diese\rond N^* \nonumber \\
\lie{\alpha}{\beta}_{\pi_N}=\lie{N^*\alpha}{\beta}_{\pi}
+\lie{\alpha}{N^*\beta}_{\pi}-N^*\lie{\alpha}{\beta}_{\pi}
\label{eq:newnir}
\end{gather*}
where $\pi_N$ is the bivector field on $X$ defined
by the relation $\pi_N\diese=\pi\diese\rond N^*$ and for any
bivector field $\ppi$ on $M$,
\begin{equation} \label{eq:pi}
\lie{\alpha}{\beta}_{\ppi}:=\derlie{\ppi\diese\alpha}(\beta)-
\derlie{\ppi\diese\beta}(\alpha)-d\big(\ppi(\alpha,\beta)\big),
\qquad \forall \alpha,\beta\in\df^1(M).
\end{equation}

\begin{prop}\label{watchingadvd}
Let $X$ be a complex manifold with associated  almost complex structure $J$.
Then  $\pi=\pire+i\piim$, where $\pire,\piim\in\sections{\wedge^2 TX}$,
is a holomorphic Poisson structure on $X$ iff
the pair $(\piim,J)$ is a Poisson Nijenhuis structure and
$\pire\diese=\piim\diese\rond J^*$.
\end{prop}

\begin{proof}
First observe that, for all $f\in\cinf{X,\CC}$ and $\alpha,\beta\in\Omega^1_\CC(X)$,
one has \[ \lie{\alpha}{f\beta}_{\pire}=\big(\pire\diese\alpha\big)(f)\beta
+f\lie{\alpha}{\beta}_{\pire} \]
and
\begin{multline*}
\lie{J^*\alpha}{f\beta}_{\piim}+\lie{\alpha}{J^*f\beta}_{\piim}
-J^*\lie{\alpha}{f\beta}_{\piim} \\
= \big(\piim\diese J^*\alpha\big)(f)\beta
+f\big(\lie{J^*\alpha}{\beta}_{\piim}+\lie{\alpha}{J^*\beta}_{\piim}-J^*\lie{\alpha}{\beta}_{\piim}\big) .
\end{multline*}
Therefore, since by Lemma~\ref{lem:2.5}, we have
 $\pire\diese=\piim\diese\rond J^*$, it suffices to check
the compatibility condition \eqref{eq:newnir}
 for $\alpha=df$ or $d\cc{f}$ and
$\beta=dg$ or $d\cc{g}$ with $f$ and $g\in\hf{X} (U)$.

An easy but cumbersome computation, making use of the relations of
Lemma~\ref{lem:1.3} and the well-known equivalences
\[ f\in\hf{X} (U) \quad \Longleftrightarrow \quad
J^* df=i df \quad \Longleftrightarrow \quad J^* d\cc{f}=-i d\cc{f}, \]
shows that the Poisson Nijenhuis compatibility of $\piim$ and $J$ is equivalent
to the closure of $\hf{X} (U)$ under the Poisson bracket of functions associated to $\pi$.

For example, $\forall f,g\in\hf{X} (U)$:
\begin{align*}
& \lie{J^* df}{dg}_{\piim}+\lie{df}{J^* dg}_{\piim}
-J^*\lie{df}{dg}_{\piim}-\lie{df}{dg}_{\pire} \\
=& \lie{i\;df}{dg}_{\piim}+\lie{df}{i\;dg}_{\piim}
-J^*d\pbim{f}{g}-d\pbre{f}{g} \\
=& 2i d\pbim{f}{g}-J^*d\pbim{f}{g}-d\pbre{f}{g} \\
=& d\pb{f}{g}+\tfrac{i}{2}J^*d\pb{f}{g}-\thalf d\pb{f}{g} \\
=& \tfrac{i}{2}\big( J^*d\pb{f}{g}-i\;d\pb{f}{g} \big) .
\qedhere \end{align*}
\end{proof}

\begin{thm}\label{PNGC} Given a complex manifold $X$ with associated
almost complex structure $J$, the following are equivalent:
\begin{enumerate}
\item\label{pq1} $\pi=\pire+i\piim\in\sections{\wedge^2T^{1,0} X}$
is a holomorphic Poisson bivector field;
\item\label{pq2} $(\piim,J)$ is a Poisson Nijenhuis structure on $X$ and
$\pire\diese=\piim\diese\rond J^*$;
\item\label{pq3} the endomorphism
\[ \JJ_\pi = \begin{pmatrix} J & \piim\diese \\ 0 & -J^* \end{pmatrix} \]
of $\TTM$ is a generalized complex structure and
$\pire\diese=\piim\diese\rond J^*$.
\end{enumerate}
\end{thm}

\begin{proof}
\fboxr{\eqref{pq1}$\iff$\eqref{pq2}} This is Proposition~\ref{watchingadvd}.
\fboxr{\eqref{pq2}$\iff$\eqref{pq3}} The equivalence follows from
 \cite[Theorem~7.6]{MR2276462} (see also \cite{Crainic2007}).
\end{proof}

\begin{rmk}
It is well known that a holomorphic Poisson structure gives
rise to a generalized complex structure. 
The holomorphic Poisson tensor is a strong Hamiltonian operator in
the sense of Liu-Weinstein-Xu \cite{MR1472888}, which deforms
the Dirac structure on $T_\CC X\oplus T_\CC^* X$
 associated to the usual complex structure seen as a generalized
complex structure \cite{Block2007,Gualtieri2007a,
MR2055289}.
\end{rmk}

\begin{rmk}
The equivalence of (a) and (b) in Theorem~\ref{PNGC} was also known to Lu \cite{private}.
\end{rmk}

\subsection{Holomorphic symplectic manifolds}

Let $(X,\omega)$ be a holomorphic symplectic manifold, where
$\omega\in\Omega^{2,0}(X)$ is the holomorphic symplectic
2-form whose corresponding holomorphic Poisson bivector field
is denoted by $\pi=\pi_R+i\pi_I\in\sections{\wedge^2 T^{1,0}X}$.
Let $\omega_R,\omega_I\in\Omega^2(X)$ be the real and imaginary
parts of $\omega$, i.e. $\omega=\omega_R+i\omega_I$.
By the holomorphic Darboux theorem both $\omega_R$ and $\omega_I$ are symplectic 2-forms.

\begin{prop}
\label{prop:2.9}
The Poisson bivector fields corresponding to
$\omega_R$ and $\omega_I$ are $4\pi_R$ and $-4\pi_I$,
respectively.
\end{prop}

\begin{proof}
The holomorphic  Darboux theorem asserts that, in a neighborhood of each point,
there exist complex symplectic coordinates $(z_1,\dots,z_n,z_1',\dots,z_n')$
so that $\omega$ can be written as \[ \omega=\sum_{k=1}^n dz_k \wedge dz_{k}' .\]
In terms of real coordinates, defined by
\[ z_k = x_k + i y_k,  \qquad z_k' = x_k' + i y_k' \]
for $k=1,\dots, n$, we have
\begin{equation}\label{eq:dar1}
\left\{ \begin{aligned} \omega_R &= \sum_{k=1}^n (dx_k \wedge dx_{k}' - dy_k \wedge dy_{k}') \\
\omega_I &= \sum_{k=1}^n (dx_k \wedge dy_{k}' + dy_k \wedge dx_{k}') \end{aligned} \right.
\end{equation}
By $\omega_R^{-1}$ and $\omega_I^{-1}$, we denote the
Poisson bivector fields corresponding to $\omega_R$ and $\omega_I$, respectively.
Thus we have
\[ \left\{ \begin{aligned}
\pi &= -\sum_{k=1}^n \frac{\partial}{\partial z_{k}}
\wedge \frac{\partial}{\partial z_{k}'} \\
\omega_R^{-1} &= -\sum_{k=1}^n \big(\frac{\partial}{\partial
x_{k}} \wedge \frac{\partial}{\partial x_{k}'}
-\frac{\partial}{\partial y_{k}} \wedge \frac{\partial}{\partial y_{k}'}\big) \\
\omega_I^{-1} &= -\sum_{k=1}^n \big(\frac{\partial}{\partial
x_{k}} \wedge \frac{\partial}{\partial y_{k}'} +
\frac{\partial}{\partial y_{k}} \wedge \frac{\partial}{\partial
x_{k}'}\big)
\end{aligned} \right. \]
On the other hand, using the relations
$\frac{\partial}{\partial z_{k}}=\frac{1}{2}(\frac{\partial}{\partial x_{k}}
-i\frac{\partial}{\partial y_{k}})$ and
$\frac{\partial}{\partial z'_{k}}=\frac{1}{2}(\frac{\partial}{\partial x'_{k}}
-i\frac{\partial}{\partial y'_{k}})$, it is simple to see that
the real and imaginary parts of  $\pi$ are given by
\[ \left\{ \begin{aligned}
\pi_{R} &= -\frac{1}{4} \sum_{k=1}^n \big( \frac{\partial}{\partial x_{k}} \wedge
\frac{\partial}{\partial x_{k}'} - \frac{\partial}{\partial y_{k}} \wedge
\frac{\partial}{\partial y_{k}'} \big) \\
\pi_{I} &= \frac{1}{4} \sum_{k=1}^n \big( \frac{\partial}{\partial x_{k}} \wedge
\frac{\partial}{\partial y_{k}'} + \frac{\partial}{\partial y_{k}} \wedge
\frac{\partial}{\partial x_{k}'} \big)
\end{aligned} \right. \]
The conclusion thus follows immediately.
\end{proof}

\subsection{Symplectic foliation}

\begin{prop}
Let $(X,\pi)$ be a holomorphic Poisson manifold, and $\pi_R$ and
$\pi_I$ the real and imaginary parts of $\pi$. Then the
symplectic foliations of $\pi_R$ and $\pi_I$ coincide,
and their leaves are exactly the holomorphic symplectic leaves of $\pi$.
\end{prop}

\begin{proof}
The relation $ \pire\diese = \piim\diese\rond J^*$ implies that $\pire$
and $\piim$ have the same symplectic leaves, for the distributions
$\pire\diese(T^*X)$ and $\piim\diese(T^*X)$ coincide.

The relation $\pi = \pire+ i \piim = \pire+  iJ \pire $ implies that,
 for all $\alpha \in (T^{0,1} X)^* $,
\[ \pi\diese(\alpha) = (\idn+iJ)\pire\diese(\alpha) .\]
 Since $ \alpha = \Re (\alpha ) + iJ^* \Re (\alpha)$, we obtain:
\[ \pi\diese(\alpha) = 2(\idn+iJ)\pire\diese\big(\Re(\alpha)\big) .\]
Taking the real part, we obtain
\[ \Re\big(\pi\diese(\alpha)\big) = 2\pi_R\diese\big(\Re(\alpha)\big) .\]
In particular, the map $T^{0,1} X\to TX$ sending an element
to its real part is an isomorphism from the distribution $\pi\diese((T^{0,1} X)^*)$ 
to the distribution $\pi_R\diese(T^* X)$, so
that the leaves associated to these distributions coincide.
\end{proof}

\section{Holomorphic Lie algebroids}

\subsection{Definition}

Holomorphic Lie algebroids were studied for various purposes in 
the literature. See \cite{MR2180064, MR1726784, MR1696093,MR1718638,MR2285039} and references
 cited there for details. Let us recall its definition below.

The tangent bundle $TX\to X$ of a complex manifold $X$ is naturally
a holomorphic vector bundle. We will denote its sheaf of holomorphic sections,
i.e. the sheaf of holomorphic vector fields, by $\Theta_X$.

Given a holomorphic vector bundle $A\xrightarrow{p}X$,
the sheaf of holomorphic sections $\shs{A}$ of $A\to X$ is the contravariant
functor which associates to an open subset $U$ of $X$ the space $\shs{A}(U)$
of holomorphic sections of $A\to X$ over $U$. Similarly, the sheaf of
 smooth sections $\shs{A}_\infty$ is the contravariant functor
 $U\to\sections{A_U}$. Clearly, $\shs{A}$ is a sheaf of $\hf{X}$-modules
while $\shs{A}_\infty$ is a sheaf of $\cinf{X}$-modules.
Moreover $\shs{A}$ is a subsheaf of $\shs{A}_\infty$.

\begin{defn}\label{def:HoloLieAlg}
A holomorphic Lie algebroid is a holomorphic vector bundle $A\to X$,
equipped with a holomorphic bundle map $A\xrightarrow{\rho}TX$,
called the anchor map, and a structure of sheaf of complex Lie algebras on $\shs{A}$,
such that
\begin{enumerate}
\item the anchor map $\rho$ induces a homomorphism of sheaves
of complex Lie algebras from $\shs{A}$ to $\Theta_X$;
\item and the Leibniz identity
\[ \lie{V}{fW}=\big(\rho(V)f\big) W+f\lie{V}{W} \]
holds for all $V,W\in\shs{A}(U)$,
 $f\in\hf{X}(U)$ and  all open subsets $U$ of $X$.
\end{enumerate}
\end{defn}

\begin{rmk}
Note that in the definition above, the last axiom implies
that the anchor map $\rho$ is automatically a holomorphic bundle
map once we assume that it is a complex bundle map.
\end{rmk}

\subsection{Underlying real Lie algebroid}

By forgetting the complex structure,
a holomorphic vector bundle $A \to X $ becomes a real (smooth) vector bundle,
and a holomorphic vector bundle map $\rho: A \to TX $
becomes a real (smooth) vector bundle map.

Let $A\to X$ be a holomorphic vector bundle whose underlying real vector
bundle is endowed with a Lie algebroid structure $(A,\rho,\lie{\cdot}{\cdot})$
such that, for any open subset $U\subset X$,
\begin{enumerate}
\item $\lie{\shs{A}(U)}{\shs{A}(U)}\subset\shs{A}(U)$
\item and the restriction of the Lie bracket
$\lie{\cdot}{\cdot}$ to $\shs{A}(U)$ is $\CC$-linear.
\end{enumerate}
Then the restriction of $\lie{\cdot}{\cdot}$ and $\rho$ 
from $\sections{A}$ to $\shs{A}$ makes $A$ a holomorphic Lie algebroid.

The following proposition states that any holomorphic Lie algebroid
can be obtained out of such a real Lie algebroid, in a unique way.

\begin{prop}
\label{prop:extension}
Given a structure of holomorphic Lie algebroid on the holomorphic
vector bundle $A\to X$ with anchor map $A\xrightarrow{\rho}TX$,
there exists a unique structure of real smooth Lie algebroid on
the vector bundle $A\to X$ with respect to the same anchor map
$\rho$ such that the inclusion of sheaves $\shs{A}\subset\shs{A}_\infty$
is a morphism of sheaves of real Lie algebras.
\end{prop}

\begin{proof}
\emph{(i) Unicity.} We first prove the unicity. Assume there exist
two such Lie algebroid structures on the vector bundle $A\to X$.
The two anchor maps would be equal. 
And for each open subset $U$ of $X$,
the two brackets would coincide on the
 subspace $\shs{A}(U)$ of $\sections{A_U}$.
Thus the two brackets would also coincide on the $\cinf{U,\RR}$-span of
$\shs{A}(U)$ inside $\sections{A_U}$. But for all trivializing open subsets
$U$ of $X$, $\shs{A}(U)$ generates $\sections{A_U}$. Hence the two
Lie algebroid structures must coincide.

\emph{(ii) Existence.} We first prove the existence of such a structure of real Lie
algebroid. Denote by $j: A \to A $ the  bundle map defining
the fiberwise complex structure on $A$.

Recall that, given a Lie algebra $\mathfrak{h}$ and a commutative algebra $F$ over a field $k$
together with a Lie algebra homomorphism $\rho:\mathfrak{h}\to\DER(F)$,
the tensor product $F\otimes_{k}\mathfrak{h}$
 is endowed with a natural Lie algebra structure over the field $k$
given by
\beq{eq:bigbra} \lie{f\otimes a}{g\otimes b}=fg\otimes\lie{a}{b}+f\rho(a)g\otimes b
-g \rho(b)f\otimes a ,\eeq
for all $a,b\in\mathfrak{h}$, $f,g\in F$.
Choose an open subset $U$ of $X$. Applying the previous general fact
to the Lie algebra $\shs{A}(U)$, the commutative algebra
$\cinf{U,\CC}$ and the anchor map $\rho:\shs{A}(U)\to\DER(\cinf{U,\CC})$,
one obtains a Lie algebra structure on $\cinf{U,\CC}\otimes_{\RR}\shs{A}(U)$.
Note that this is a real Lie algebra, since $\shs{A}(U)\to\DER(\cinf{U,\CC})$
is $\RR$-linear but \emph{not} $\CC$-linear.

Now, for any holomorphic function $h\in\hf{X}(U)$, it follows from 
Eq.~\eqref{eq:bigbra} and the Leibniz identity for the holomorphic
Lie algebroid $A\to X$ that
\begin{align*}
\lie{f\otimes a}{gh\otimes b-g\otimes hb} =& \;
fgh\otimes \lie{a}{b} + fg \rho(a)h\otimes b + f \big( \rho(a)g \big) h\otimes b \\
& -fg\otimes h\lie{a}{b}- fg\otimes \big(\rho(a)h\big)b - f\big(\rho(a)g\big)\otimes hb .
\end{align*}

In other words, the elements of type $fh\otimes a-f\otimes ha$, with
$f\in\cinf{U,\CC}$, $a\in\shs{A}(U)$ and $h\in\hf{X}(U)$,
generate an ideal of the Lie algebra $\cinf{U,\CC}\otimes_{\RR}\shs{A}(U)$.
As a consequence, the Lie bracket given by Eq.~\eqref{eq:bigbra} induces a
Lie algebra structure (over $\RR$) on the quotient
 $\cinf{U,\CC}\otimes_{\hf{X}(U)}\shs{A}(U)$
 of $\cinf{U,\CC}\otimes_{\RR}\shs{A}(U)$ by the aforementioned ideal.

There is a natural map $\cinf{U,\CC}\otimes_{\hf{X}(U)}\shs{A}(U)
\injection\sections{A_U}$ mapping $f\otimes a$ to $\Re(f)\, a+ \Im(f)\, j(a)$,
for all $f\in\cinf{U,\CC}$ and $a\in\shs{A}(U)$. This is actually an
isomorphism if the open subset $U$ of $X$ is trivializing the
holomorphic bundle $A|_U\to U$. Indeed, any smooth section of $A\to X$
over $U$ can be written as a linear combination $\sum_{k=1}^{m}f_k
a_k$, where $a_k\in\shs{A}(U)$ and $f_k$ is a smooth $\CC$-valued
function on $U$. Therefore, $\sections{A_U}$ is a Lie algebra.
By construction, the previous Lie
bracket restricts to a $\CC$-linear bracket on $\shs{A}(U)$ and
is a Lie-Rinehart algebra over $\cinf{U}$. Hence, $A|_U$ is a smooth 
Lie algebroid. By the unicity in (i), one obtains a 
smooth Lie algebroid $A\to X$ by gluing $A|_U$ together using an open
covering $\{U_i\}$ of $X$. 
\end{proof}

In the sequel, we will use $A_R$ to denote the
underlying real Lie algebroid of a holomorphic Lie algebroid $A$.
 When referring to holomorphic Lie algebroids,
we either use Definition~\ref{def:HoloLieAlg},
or the equivalent one as in Proposition~\ref{prop:extension}.
In particular, by saying that
a real Lie algebroid is a holomorphic Lie algebroid,
we mean that it is a holomorphic vector bundle and
its Lie bracket on smooth sections restricts to a $\CC$-linear
bracket on $\shs{A}(U)$, for all open subset $U\subset X$.

\subsection{Underlying imaginary Lie algebroid}
\label{sec:4.1}

Assume that $(A\to X,\rho,\lie{\cdot}{\cdot})$ is
a holomorphic Lie algebroid. Consider
the bundle map $j: A \to A$ defining
the fiberwise complex structure on $A$.
We compute the Nijenhuis torsion of $j$ by considering $A$ as
a real Lie algebroid $A_R$.

\begin{prop}\label{j-Nijenhuis}
Let $(A\to X,\rho,\lie{\cdot}{\cdot})$ be a holomorphic Lie algebroid
and $j:A_R\to A_R$ its associated endomorphism.
Then the Nijenhuis torsion of $j$ vanishes.
\end{prop}

\begin{proof}
Since $T(j)$ is a section in $\gm (\wedge^2 A_R^* \otimes A_R)$, 
it suffices to prove that $T(j)(V,W)=0$ for any holomorphic
sections $V,W\in\shs{A}(U)$, where $U\subset X$ is
any open subset. This however follows
immediately from the $\CC$-linearity of $\lie{\cdot}{\cdot}$ on $\shs{A} (U)$:
\beq{zzz} \lie{V}{jW}=j\lie{V}{W}=\lie{jV}{W}, \qquad\forall V,W\in\shs{A}(U) .\eeq
\end{proof}

Since the Nijenhuis torsion of $j:A_R\to A_R$ vanishes,
one can define a new (real) Lie algebroid structure on
$A$, denoted by $(A\to X,\rho_{j},\lie{\cdot}{\cdot}_j)$,
where the anchor $\rho_{j}$ is $\rho\rond j$ and
the bracket on $\sections{A}$ is
given by \cite{MR1421686}
\[ \lie{V}{W}_{j}=\lie{jV}{W}+\lie{V}{jW}-j\lie{V}{W}, \qquad\forall V,W\in\sections{A} .\]

In the sequel, $(A\to X,\rho_{j},\lie{\cdot}{\cdot}_j)$
will be called the \emph{underlying imaginary Lie algebroid}
and denoted by $A_I$. It is known that 
\begin{equation}
j: A_I\to A_R
\end{equation}
 is a Lie algebroid isomorphism \cite{MR1421686}.

\subsection{Associated complex Lie algebroids}

Let $(A\to X,\lie{\cdot}{\cdot},\rho)$ be a holomorphic Lie algebroid.
Complexifying its underlying real Lie algebroid (which was described in
Proposition~\ref{prop:extension}) by extending the anchor map and the
Lie bracket $\CC$-linearly,
we obtain a complex Lie algebroid:
\[ \xymatrix{ A\cmplx \ar[d] \ar[r]^{\canchor} & T\cmplx X \ar[dl] \\ X & } ,\]
where $A_\CC =A\otimes \CC$.
Extending  $\CC$-linearly the bundle map $j:A\to A$, 
we obtain a  map of complex vector bundles
$j: A\cmplx\to A\cmplx$ such that $j^2=-\idn$.
Let $A^{1,0}\to X$ and $A^{0,1}\to X$
be its eigenbundles with eigenvalues $i$ and $-i$ respectively.
It follows from Proposition~\ref{j-Nijenhuis} that
$\sections{A^{1,0}}$ and $\sections{A^{0,1}}$ are Lie subalgebras
of $\sections{A\cmplx}$.
Hence $A^{1,0}$ and $A^{0,1}$ are complex Lie subalgebroids of $A\cmplx$.
Note that the map
\begin{equation}
A \to A^{1,0} : a \mapsto \tfrac{1}{2}(a - i j (a))
\end{equation}
is an isomorphism of complex vector bundles.
Hence by pulling back all the structures,
one obtains a complex Lie algebroid structure on
the same complex vector bundle $A\to X$.
Similarly, $a \to \frac{1}{2}(a + i j (a))$
is an isomorphism of complex vector bundles
$A\to A^{0,1}$. Hence $A\to X$ admits
another complex Lie algebroid structure.
 The following proposition describes
these complex Lie algebroids explicitly.
Its proof is a simple computation and is left
to the reader.

\begin{prop}
Given a holomorphic Lie algebroid $(A\to X,\lie{\cdot}{\cdot},\rho)$,
let $(A\to X,\lie{\cdot}{\cdot}_{1,0},\rho_{1,0})$
and $(A\to X,\lie{\cdot}{\cdot}_{0,1},\rho_{0,1})$
be its associated complex Lie algebroids corresponding
to $A^{1,0}$ and $A^{0,1}$, respectively.
Then
\begin{align}
\rho_{1,0}&=\tfrac{1}{2}(\rho+i\rho_j),  & \lie{\cdot}{\cdot}_{1,0}&=\tfrac{1}{2}(\lie{\cdot}{\cdot}+j\lie{\cdot}{\cdot}_j) \\
\rho_{0,1}&=\tfrac{1}{2}(\rho-i\rho_j),  & \lie{\cdot}{\cdot}_{0,1}&=\tfrac{1}{2}(\lie{\cdot}{\cdot}-j\lie{\cdot}{\cdot}_j) .
\end{align}
\end{prop}

%\begin{proof}
%By construction, the anchor map $\rho_{1, 0}$ is given by
% \[ \rho_{1, 0}(a) = \rho (  \frac{1}{2}(a + i j (a)) ), \ \ \ \forall a\in A .\]
%Hence, it follows that $\rho_{1, 0}=\frac{1}{2}(\rho + i \rho_j )$.
%
%As for  the bracket $ [\cdot, \cdot ]_{1, 0}$,
% for any sections $V,W \in  \sections{A}$:
%  \[  \frac{1}{2}( [V,W ]_{1, 0} + i j ([V,W ]_{1, 0} )) =
%     [  \frac{1}{2}(V + i j V) , \frac{1}{2} (W + i j W) ] ,\]
%where the bracket on the right hand side is the $\CC$-linear  bracket on
%$\sections{A_\CC}$. Hence we have
% \[ \frac{1}{2}( [V,W ]_{1, 0} + i j ([V,W ]_{1, 0} ))=\frac{1}{4}([V,W]- [jV,jW]+ i ( [jV,W ]+[V,jW ]) ) .\]
%Taking the real parts,  one obtains
% \[ [V,W ]_{1, 0} = \frac{1}{2} ([V,W]+ j [V,W]_j) .\]
%The other relations can be proved in a similar fashion.
%\end{proof}

\subsection{Cotangent bundle Lie algebroids of holomorphic Poisson manifolds}

In this section, as an example  we consider the cotangent  bundle Lie algebroid
of  a holomorphic Poisson manifold and identify various Lie algebroids
associated to it, which  will be used later on.

Assume that $(X, \pi)$ is a holomorphic Poisson manifold, where
$\pi=\pi_R+i \pi_I\in \gm (\wedge^2 T^{1, 0}X)$. Let
$A=(T^*X)_\pi$ be its corresponding holomorphic Lie algebroid,
which can be  defined in a similar way as in the smooth case. 
To be more precise, let $\Phi$ and $\Psi$, respectively, be the
holomorphic bundle maps
\[ \Phi: TX \to T^{1, 0}X, \ \Phi=  \frac{1}{2} (1-iJ)  \]
and
\[ \Psi: T^*X\to (T^{1, 0}X)^*, \Psi=1-iJ^* ,\]
where $J$ is the almost complex structure on $X$.
Define the anchor $\rho: (T^*X)_\pi\to TX$ to be
$\rho=\Phi^{-1}\rond \pi^{\#} \rond \Psi$ and the bracket
\[ [\alpha, \beta]_\pi =  L_{\rho \alpha } \beta - L_{\rho \beta }
   \alpha  -   {\rm d} (\rho \alpha,   \beta ) \]
 $\forall \alpha,\beta  \in \Gamma ( T^*X|_U)$ holomorphic.

Its associated complex Lie algebroid $A^{1,0}$ will be denoted $(T^{1,0}
X)^*_\pi$ since its underlying complex vector bundle is $(T^{1,0}
X)^*$. The following lemma is obvious.

\begin{lem}
For the associated complex Lie algebroid $(T^{1,0} X)^*_\pi$, 
the anchor map is 
\[ (T^{1,0} X)^*\xrightarrow{\pi\diese}T\cmplx X \]
 and the bracket on $\Omega^{1, 0}(X)$ is given by
\[ [\xii, \etaa]=L_{\pi\diese\xii}\etaa-L_{\pi\diese\etaa}\xii
-\partial \pi (\xii, \etaa ), \qquad \forall \xii , \etaa \in
\Omega^{1, 0}(X) .\]
\end{lem}

The following proposition describes the underlying real and
imaginary Lie algebroids of $(T^*X)_\pi$.

\begin{prop}
\label{prop:3.7}
Let $(X, \pi )$ be a holomorphic Poisson manifold,
 where
$\pi=\pi_R+i \pi_I\in \gm (\wedge^2 T^{1, 0}X)$. Then
the underlying real and
imaginary Lie algebroids of $(T^*X)_\pi$ are isomorphic
to $(T^*X)_{4\pi_R}$ and $(T^*X)_{4\pi_I}$, respectively.
\end{prop}
\begin{proof}
First, let us consider the anchor map $\rho :(T^* X)_\pi\to TX$
as a bundle map of real vector bundles.
Clearly we have
\[ \rho=\Phi^{-1}\rond \pi^{\#}\rond \Psi
=\Phi^{-1} \rond \pi^\#_R \rond ( 1 -  i J^* )\rond (1-iJ^*) =4 \pi_R^{\#} .\]

Now we consider the bracket. For this purpose,  let $A_R$ be  its
underlying real Lie algebroid and $A_\CC$ the complexification
of $A_R$.
Note that for any holomorphic functions $f, g \in \hf{X}(U)$,
we have
\[ [df, dg]_\pi=d\{f, g\} ,\]
where both sides are considered as sections of $A_\CC$.
It thus follows, from Corollary \ref{cor:2.4},
that
\[ \Re [df, dg]_\pi=d\Re \{f, g\}=4d\{\Re  f,  \Re g\}_R=
4[{\rm d} \Re (f) , {\rm d} \Re (g) ]_{\pi_R} .\]
Therefore, $A_R$ is isomorphic to $(T^*X)_{4\pi_R}$.

Finally note that the Nijenhuis structure $j: A_R \to A_R$
in Section \ref{sec:4.1} is simply $J^*$. It thus follows that
$(A_R)_j$ is isomorphic to $(T^*X)_{4\pi_I}$ using Theorem~\ref{PNGC}.
\end{proof}

\subsection{Holomorphic Lie-Poisson structures}

Similar to the smooth case,
there is also another equivalent definition of holomorphic Lie
algebroids. The proof is similar to the smooth case, and is
left to the reader.
Note that the complex dual $\CA{A}$ of a holomorphic vector bundle
$A \to X $ is again a complex manifold, which is
also a holomorphic vector bundle over $X$. 
We denote by $p:\CA{A}\to X$ the projection onto
the base manifold. There is a one-one correspondence between holomorphic
sections $V\in\shs{A}(U)$ and fiberwise-linear holomorphic
functions $l_{V}$ on $\CA{A|_U}$: $\forall\alpha\in\CA{A|_U}$
\[ l_{V} (\alpha) = \alpha (\left. V \right|_{p(\alpha)} ) .\]

\begin{prop}
\label{prop:dualPoisson}
Let $A \to X $ be a holomorphic vector bundle.
The following are equivalent:
\begin{enumerate}
\item $A$ is a holomorphic Lie algebroid;
\item there exists a fiberwise-linear holomorphic Poisson structure on $\CA{A}$.
\end{enumerate}
Here the Lie algebroid structure on $(A, \rho, [\cdot, \cdot])$
and the Poisson structure on $\CA{A}$ are related by
the following equations:
\[ \begin{array}{rcl} \{p^* f, l_{V}\} & = & p^* \big( \rho(V)(f)\big)  \\
  \{l_{V}, l_{W} \}& = & l_{[V,W]}
 \end{array} \]
for any $ V,W \in \shs{A} (U) $ and $f \in {\mathcal O}_X (U)$,
where $ p: \CA{A} \to X$ is the projection.
\end{prop}

Summarizing the discussions above, we get the main result
of this section.

\begin{thm}\label{thm:dualPoisson}
Let $A \to X$ be holomorphic vector bundle.
There is a one-to-one correspondence between:
\begin{enumerate}
\item holomorphic Lie algebroid structures on $A \to X$,
\item fiberwise-linear holomorphic Poisson structures on $\CA{A}$,
\end{enumerate}
\end{thm}

\subsection{Real and imaginary parts of a holomorphic Lie-Poisson structure}

Any complex vector space $V$ can be equivalently thought of as a real vector space
with an $\RR$-linear endomorphism $j$, such that $j^2=-1$, representing the multiplication
by the imaginary number $\sqrt{-1}$.

Given a complex vector space $V$, its complex dual space is the set of morphisms
$\Hom_{\CC}(V,\CC)$ from $V$ to $\CC$ in the category of complex vector spaces.
Similarly, its real dual space is the set of morphisms $\Hom_{\RR}(V,\RR)$
from $V$ to $\RR$ in the category of real vector spaces.
Clearly, $\Hom_{\CC}(V,\CC)$ is a vector space over $\CC$
while $\Hom_{\RR}(V,\RR)$ is a vector space over $\RR$.

The map
\[ \Hom_{\RR}(V,\RR) \to \Hom_{\CC}(V,\CC) : f \mapsto (1-ij^*)f \]
is an isomorphism of \emph{real} vector spaces.
Indeed, $g\in\Hom_{\CC}(V,\CC)$ if and only if $g=(1-ij^*)f$ with $f(=\Re\rond g)\in\Hom_{\RR}(V,\RR)$.

%\begin{align*}
%& g\in\Hom_{\CC}(V,\CC) \\
%\Leftrightarrow \qquad & g\in\Hom_{\RR}(V,\CC) \text{ and } g(jv)=ig(v),\;\forall v\in V \\
%\Leftrightarrow \qquad & g\in\Hom_{\RR}(V,\CC) \text{ and } \Re\rond g\rond j=-\Im\rond g \\
%\Leftrightarrow \qquad & g=(1-ij^*)f \text{ with } f(=\Re\rond g)\in\Hom_{\RR}(V,\RR) .
%\end{align*}

Given a complex vector bundle $A\to X$,
 we denote its complex and real dual bundles
by $\Hom_{\CC}(A,\CC)\to X$ and $\Hom_{\RR}(A,\RR)\to X$,
respectively. Applying the previous isomorphism fiberwise yields
the isomorphism of real vector bundles $\Psi=1-ij^*$:
\beq{eq:psiq} \xymatrix{
\Hom_{\RR}(A,\RR) \ar[d]_{p_{\RR}}\ar[r]^{\Psi} & \Hom_{\CC}(E,\CC) \ar[d]^{p_{\CC}} \\
X \ar[r]_{\idn} & X .
} \eeq
Here $p_{\CC}$ and $p_{\RR}$ denote
the projections of the vector bundles $\CA{A},\RA{A}$
onto their base $X$. Note that $\Psi\inv(\xi)=\Re\rond\xi$.

We consider a holomorphic Lie algebroid $(A\to X,\rho,\lie{\cdot}{\cdot})$.
According to Proposition~\ref{prop:dualPoisson}, the complex dual
bundle $\CA{A}$ is a fiberwise linear holomorphic Poisson manifold,
whose holomorphic Poisson tensor is denoted by $\pi$.
Let $\pire$ and $\piim$ be its real and imaginary parts.
Then $\picre:=\Psi\inv_*\pire$ and $\picim:=\Psi\inv_*\piim$
are fiberwise $\RR$-linear Poisson tensors on
the real dual bundle $\RA{A}$.
These Poisson structures therefore correspond to real Lie
algebroids on $A\to X$, which are denoted by
$(A\to X,\rho_{\Re},\lie{\cdot}{\cdot}_\Re)$
or $A_\Re$ 
and $(A\to X,\rho_{\Im},\lie{\cdot}{\cdot}_\Im)$ 
or $A_\Im$, respectively.
The following Proposition identifies these Lie algebroids
with those discussed in Section~\ref{sec:4.1}.

\begin{prop}\label{prop:realparts}
Let $(A\to X,\rho,\lie{\cdot}{\cdot})$ be a holomorphic Lie algebroid.
\begin{enumerate}
\item The Lie algebroid $(A\to X,4\rho_{\Re},4\lie{\cdot}{\cdot}_\Re )$ is
isomorphic to the real Lie algebroid $A_R$;
% $(A\to X,\rho,\lie{\cdot}{\cdot})$.
\item The Lie algebroid $(A \to X,-4\rho_{\Im},-4\lie{\cdot}{\cdot}_{\Im})$
is isomorphic to the imaginary Lie algebroid $A_I$.
% $(A\to X,\rho_{j},\lie{\cdot}{\cdot}_{j})$.
\end{enumerate}
\end{prop}

\begin{proof}
First, we fix some notations.
Any section $V\in\sections{A}$ can be thought of as a fiberwise $\CC$-linear
(resp. $\RR$-linear) function on $\CA{A}$ (resp. $\RA{A}$), which we denote by
$l_V $ (resp. $l_V'$).

For all $V\in\sections{A}$ and $f\in\RA{A}$, we have
\[ \big(\Psi^*l_V\big)(f)=l_V(f-i(j^*f))=f(V)-i\;f(jV)=l'_V(f)-i\;l'_{jV}(f) .\]
Thus
\beq{eq:isom}
\left\{\begin{aligned}
\Re(\Psi^*l_V)&=l_V' \\ \Im(\Psi^*l_V)&=-l_{jV}'
.\end{aligned}\right.
\eeq

First we look at the anchor map $\rho_{\Re}$.
Given any $x\in X$, $\alpha\in A_x$ and $\eta\in T^*_x X$,
to compute $\pairing{\eta}{\rho_{\Re}(\alpha)}$,
we choose an open neighborhood $U\subset X$ containing $x$,
a holomorphic section $V\in\shs{A}(U)$ through $\alpha$,
and a holomorphic function $f\in\hf{X}(U)$
with $\diff\big(\Re(f)\big)|_x =\eta$.

Consider the relation
\[ p_{\CC}^*\big(\rho(V)f\big)=\pb{p_{\CC}^* f}{l_V} \]
from Proposition~\ref{prop:dualPoisson}.
Taking the real part, we obtain, by Corollary~\ref{cor:2.4}:
\[ p_{\CC}^*\big(\rho(V)\Re(f)\big) =
4 \pbre{p_{\CC}^*\Re(f)}{\Re(l_V)} .\]
 Applying $\Psi^*$ to both
sides and using Eqs.~\eqref{eq:psiq}~and~\eqref{eq:isom}, we have
\begin{multline*}
p_{\RR}^* \big(\rho(V)\Re(f)\big)
= \Psi^* p_{\CC}^* \big(\rho(V)\Re(f)\big)
= 4 \Psi^* \pbre{p_{\CC}^*\Re(f)}{\Re(l_V)} \\
= 4 \pbcre{\Psi^* p_{\CC}^*\Re(f)}{\Re(\Psi^* l_V)}
= 4 \pbcre{p_{\RR}^*\Re(f)}{l'_V}
= 4 p_{\RR}^*\big(\rho_{\Re}(V)\Re(f)\big)
.\end{multline*}
Hence it follows that
\[ \rho(V)\Re(f)=4\rho_{\Re}(V)\Re(f) \]
and
\[ \pairing{\eta}{4\rho_{\Re}(\alpha)}=\pairing{\eta}{\rho(\alpha)} .\]

We now identify the anchor map $\rho_\Im$.
Taking the imaginary part of the relation
\[ p_{\CC}^*\big(\rho(V)f\big)=\pb{p_{\CC}^* f}{l_V} \]
from Proposition~\ref{prop:dualPoisson}, we obtain
\[ p_{\CC}^*\big(\rho(V)\Im(f)\big)=\Im\pb{p_{\CC}^*f}{l_V}=-4\pbim{p_{\CC}^*\Im(f)}{\Im(l_V)} \]
by Corollary~\ref{cor:2.4}. Applying $\Psi^*$ to both sides and
making use of \eqref{eq:psiq} and \eqref{eq:isom}, we get:
\begin{multline*}
p_{\RR}^*\big(\rho(V)\Im{(f)}\big)=\Psi^*p_{\CC}^*\big(\rho(V)\Im(f)\big)
=-4\Psi^*\pbim{p_{\CC}^*\Im(f)}{\Im(l_V)} \\ =-4\pbcim{\Psi^*p_{\CC}^*\Im(f)}{\Psi^*(l_V)}
=4\pbcim{p_{\RR}^*\Im(f)}{l'_{jV}}
.\end{multline*}
It follows that \[ \rho(V)\Im(f)=4\rho_{\Im}(jV)\Im(f) ,\]
which implies that $4\rho_{\Im}=-\rho_j$.

 We now consider the Lie brackets.  
Choose an open neighborhood $U\subset X$ such that the holomorphic
vector bundle $A|_U\to U$ is trivial.
Since the module of smooth sections of $A|_U$ is generated (over $\cinf{U,\RR}$)
by the holomorphic sections, it suffices to show that
the bracket $4[V,W]_\Re $ (resp. $-4[V, W]_\Im$) is equal
to $\lie{V}{W}$ (resp. $\lie{V}{W}_j$), for any two
holomorphic sections $V,W\in\shs{A}(U)$.

According to Proposition~\ref{prop:dualPoisson}, we have
\beq{eq:dualPoi2} \pb{l_{V}}{l_{W}} = l_{\lie{V}{W}}
\qquad\forall V,W\in\shs{A}(U) .\eeq
By Corollary~\ref{cor:2.4}, we obtain
\[ \Re(l_{\lie{V}{W}}) = \Re\pb{l_{V}}{l_{W}} = 4\pbre{\Re(l_V)}{\Re(l_W)} .\]
Therefore, applying $\Psi^*$ to both sides, we get
\[ \Re(\Psi^*l_{\lie{V}{W}}) = 4\Psi^*\pbre{\Re(l_V)}{\Re(l_W)}
= 4\pbcre{\Re(\Psi^*l_V)}{\Re(\Psi^*l_W)} \]
and, using \eqref{eq:isom}, we obtain
\[ l'_{\lie{V}{W}}=4\pbcre{l'_V}{l'_W} .\]
Hence it follows that $4\lie{V}{W}_{\Re}=\lie{V}{W}$.

We now turn our attention to the Lie bracket $\lie{\cdot}{\cdot}_{\Im}$.
Taking the imaginary part of Eq.~\eqref{eq:dualPoi2}
and applying $\Psi^*$ to both sides, it follows from Corollary~\ref{cor:2.4} that
\[ -4\pbcim{\Im(\Psi^*l_V)}{\Im(\Psi^*l_W)} = \Im(\Psi^*l_{\lie{V}{W}}) .\]
Hence, using Eq.~\eqref{eq:isom}, we obtain
\[ 4 l'_{\lie{jV}{jW}_{\Im}} = 4 \pbcim{l'_{jV}}{l'_{jW}} = l'_{j\lie{V}{W}} .\]
Therefore, \[ 4\lie{jV}{jW}_{\Im}=j\lie{V}{W} \] and
\[ \lie{V}{W}_j = j\inv\lie{jV}{jW} = -4\lie{V}{W}_{\Im} .\]
This completes the proof.
\end{proof}

\begin{rmk}
In particular, given a Lie algebra $\mathfrak{g}$ over $\CC$,
its complex dual $\Hom_{\CC}(\mathfrak{g},\CC)$ is a holomorphic
Poisson manifold.
The isomorphism $\Hom_{\RR}(\mathfrak{g},\RR)\xrightarrow{\Psi}
\Hom_{\CC}(\mathfrak{g},\CC)$ maps the Lie-Poisson structure on
$\Hom_{\RR}(\mathfrak{g},\RR)$ corresponding to the Lie algebra
bracket $v\otimes w\mapsto \tfrac{1}{4}\lie{v}{w}$ (resp.
$v\otimes w\mapsto -\tfrac{1}{4}j\lie{v}{w}$) on $\mathfrak{g}$
to the real (resp. imaginary) part of the holomorphic Poisson
structure on $\Hom_{\CC}(\mathfrak{g},\CC)$.
Here $j:\mathfrak{g}\to\mathfrak{g}$ is the $\RR$-linear operator
representing the scalar multiplication by $\sqrt{-1}\in\CC$.

This is an immediate consequence of the relations
\[ \pbcre{l'_V}{l'_W}=l'_{\tfrac{1}{4}\lie{V}{W}} \qquad
\pbcim{l'_V}{l'_W}=l'_{-\tfrac{1}{4}\lie{V}{W}_j} \]
and the following fact: since here, the holomorphic vector bundle
$A\to X$ is reduced to the vector space $\mathfrak{g}$, all elements
$V,W$ of $\mathfrak{g}$ are automatically holomorphic sections
and, we have
\[ \lie{V}{W}_j=\lie{jV}{W}+\lie{V}{jW}-j\lie{V}{W}=j\lie{V}{W} \]
since the restriction of the Lie bracket to the holomorphic sections is $\CC$-linear.
\end{rmk}

\section{Holomorphic Lie algebroid cohomology and holomorphic Poisson cohomology}

\subsection{Matched pairs}
The notion of matched pairs of Lie algebroids was introduced
by Lu in her classification of Poisson group actions \cite{MR1430434},
and further studied by Mokri \cite{MR1460632} and Mackenzie \cite{Mackenzie2007}.

\begin{defn}
A \emph{matched pair} of Lie algebroids
is a pair of (complex or real) Lie algebroids $A$ and $B$ over
the same base manifold $M$, where $B$ is an $A$-module and
$A$ is a $B$-module such that the following identities hold:
\begin{gather}
\lie{a(X)}{b(Y)} = -a\big(\nabla_Y X\big)+b\big(\nabla_X Y\big)
, \label{first} \\
\nabla_X\lie{Y_1}{Y_2} = \lie{\nabla_X Y_1}{Y_2} +
\lie{Y_1}{\nabla_X Y_2} + \nabla_{\nabla_{Y_2} X}Y_1 -
\nabla_{\nabla_{Y_1} X} Y_2
, \label{second} \\
\nabla_Y\lie{X_1}{X_2} = \lie{\nabla_Y X_1}{X_2} +
\lie{X_1}{\nabla_Y X_2} + \nabla_{\nabla_{X_2} Y} X_1 -
\nabla_{\nabla_{X_1}Y}X_2 , \label{third}
\end{gather}
where $X_1,X_2,X\in\sections{A}$ and $Y_1,Y_2,Y\in\sections{B}$.
Here, $a$ and $b$ are the anchor maps of $A$ and $B$, respectively,
and $\nabla$ denotes both the representation
\[ \sections{A}\otimes\sections{B}\to\sections{B}:(X,Y)\mapsto\nabla_X Y \]
of $A$ on $B$ and the representation
\[ \sections{B}\otimes\sections{A}\to\sections{A}:(Y,X)\mapsto\nabla_Y X \]
of $B$ on $A$.
\end{defn}

\begin{thm}[\cite{MR1460632,Mackenzie2007}]
\label{thm:5.2}
Given a matched pair $(A,B)$ of Lie algebroids, there is a Lie
algebroid structure $A\bowtie B$ on the direct sum vector bundle
$A\oplus B$, with anchor $c(X\oplus Y)=a(X)+b(Y)$ and bracket
\beq{6bis} \lie{X_1\oplus Y_1}{X_2\oplus Y_2} = \big(
\lie{X_1}{X_2} + \nabla_{Y_1}X_2 - \nabla_{Y_2}X_1 \big) \oplus
\big( \lie{Y_1}{Y_2} + \nabla_{X_1}Y_2 - \nabla_{X_2}Y_1 \big) .
\eeq
Conversely, if $A\oplus B$ has a Lie algebroid structure for
which $A\oplus 0$ and $0\oplus B$ are Lie subalgebroids, then the
representations $\nabla$ defined by
\[ \lie{X\oplus 0}{0\oplus Y} = -\nabla_Y X\oplus\nabla_X Y \]
endow the pair $(A,B)$ with a matched pair structure.
\end{thm}
See \cite{Mackenzie2007} for more details.

\begin{ex}\label{t01t10}
Let $X$ be a complex manifold. Then $(T^{0,1} X,T^{1,0} X)$ is a
matched pair, where the actions are given by
\[ \nabla_{X^{0,1}}X^{1,0}=\pr^{1,0}\lie{X^{0,1}}{X^{1,0}} \qquad \text{ and } \qquad
\nabla_{X^{1,0}}X^{0,1}=\pr^{0,1}\lie{X^{1,0}}{X^{0,1}}, \] for all
$X^{0,1}\in\XX^{0,1}(X)$ and $X^{1,0}\in\XX^{1,0}(X)$. Hence $T^{0,1}
X\bowtie T^{1,0}X$ and $T\cmplx X$ are isomorphic as complex Lie
algebroids.

More generally, given a holomorphic Lie algebroid $A$,
the pair $(A^{0,1},A^{1,0})$
is a matched pair of Lie algebroids and $A^{0,1}\bowtie A^{1,0}$
is isomorphic, as a complex Lie algebroid, to $A\cmplx$.
\end{ex}

Let $A$ and $B$ be two (complex or real) Lie algebroids over
the same base manifold $M$. Assume $B$ is an $A$-module and
$A$ is a $B$-module, both representations being abusively denoted by the same symbol $\nabla$.
And define
\begin{gather*}
F(X;Y) := \lie{a(X)}{b(Y)} + a\big(\nabla_Y X\big) - b\big(\nabla_X Y\big), \\
S(X;Y_1,Y_2) := \lie{\nabla_X Y_1}{Y_2} + \lie{Y_1}{\nabla_X Y_2} - \nabla_X\lie{Y_1}{Y_2}
+ \nabla_{\nabla_{Y_2} X} Y_1 - \nabla_{\nabla_{Y_1} X} Y_2, \\
T(Y;X_1,Y_2) := \lie{\nabla_Y X_1}{X_2} + \lie{X_1}{\nabla_Y X_2} - \nabla_Y\lie{X_1}{X_2}
+ \nabla_{\nabla_{X_2} Y} X_1 - \nabla_{\nabla_{X_1}Y}X_2 ,
\end{gather*}
where $a$ and $b$ are the respective anchor maps of $A$ and $B$,
while $X_1,X_2,X\in\sections{A}$ and $Y_1,Y_2,Y\in\sections{B}$.

The following result can be verified directly.

\begin{lem}\label{tensoriality}
For any (complex or real-valued) function $f$ on $M$, we have:
\begin{align*}
F(fX;Y)&=fF(X;Y) & F(X;fY)&=fF(X;Y) \\
S(fX;Y_1,Y_2)&=fS(X;Y_1,Y_2) & T(fY;X_1,X_2)&=fT(Y;X_1,X_2)
\end{align*}
and
\begin{gather*}
S(X;fY_1,Y_2)=fS(X;Y_1,Y_2)+F(X;Y_2)(f)Y_1 \\
T(Y;fX_1,X_2)=fT(Y;X_1,X_2)-F(X_2;Y)(f)X_1 .
\end{gather*}
Moreover, $S$ and $T$ are skew-symmetric in their last two arguments.
\end{lem}

\subsection{Cohomology of a matched pair}

\begin{prop}\label{doublecomplex}
Let $A$ and $B$ be a pair of Lie algebroids over $M$ with mutual actions
$\nabla$. The pair $(A,B)$ is a matched pair iff the diagram \beq{eq:double}
\xymatrix{ \sections{\wedge^k A^*\otimes\wedge^l B^*}
\ar[r]^{\partial_A} \ar[d]_{\partial_B} &
\sections{\wedge^{k+1} A^*\otimes\wedge^l B^*} \ar[d]^{\partial_B} \\
\sections{\wedge^k A^*\otimes\wedge^{l+1} B^*} \ar[r]_{\partial_A} &
\sections{\wedge^{k+1} A^*\otimes\wedge^{l+1} B^*} }
\eeq
commutes, where $\partial_A$ and $\partial_B$
denote the Lie algebroid cohomology differential operators
of $A$ with values in the module $\wedge^\bullet B^*$
and of $B$ with values in the module $\wedge^\bullet A^*$, respectively.

Here, if $\alpha\in\sections{\wedge^k A^*\otimes\wedge^l B^*}$, $\partial_A$ and $\partial_B$ are given by
\beq{partialA}\begin{aligned}
& \big(\partial_A\alpha\big)(A_0,\dots,A_k,B_1,\dots,B_l) \\
= \; & \sum_{i=0}^k (-1)^i \Big(a(A_i)\alpha(A_0,\dots,\widehat{A_i},\dots,A_k,
B_1,\dots,B_l) \\
& - \sum_{j=1}^l \alpha(A_0,\dots,\widehat{A_i},\dots,A_k,
B_1,\dots, \nabla_{A_i}B_j,\dots,B_l)\Big) \\
& + \sum_{i<j} (-1)^{i+j} \alpha(\lie{A_i}{A_j},A_0,\dots,\widehat{A_i},
\dots,\widehat{A_j},\dots,A_k,B_1,\dots,B_l)
\end{aligned}\eeq
and
\beq{partialB}\begin{aligned}
& \big(\partial_B\alpha\big)(A_1,\dots,A_k,B_0,\dots,B_l) \\
= \; & \sum_{i=0}^l (-1)^i \Big(b(B_i)\alpha(A_1,\dots,A_k,
B_0,\dots,\widehat{B_i},\dots,B_l) \\
& - \sum_{j=1}^k \alpha(A_1,\dots,\nabla_{B_i}A_j,\dots,A_k,
B_0,\dots,\widehat{B_i},\dots,B_l)\Big) \\
& + \sum_{i<j} (-1)^{i+j} \alpha(A_1,\dots,A_k,
\lie{B_i}{B_j},B_0,\dots,\widehat{B_i},\dots,\widehat{B_j},\dots,B_l) .
\end{aligned}\eeq
\end{prop}

\begin{proof}
\fboxr{$\Rightarrow$}
If the pair $(A,B)$ is a matched pair, then the direct sum $A\oplus B$ is a Lie algebroid with bracket given by \eqref{6bis}. And the corresponding Lie algebroid differential
\[ \sections{\wedge^{\bullet}(A\oplus B)^*}\xrightarrow{d_{A\bowtie B}}
\sections{\wedge^{\bullet+1}(A\oplus B)^*} ,\]
defined by
\begin{multline*}
(d_{A\bowtie B}\alpha )(C_0,\dots,C_n)
=\sum_{i=0}^n(-1)^i c(C_i) \big(\alpha(C_0,\dots,\widehat{C_i},\dots,C_n)\big) \\
+\sum_{i<j}(-1)^{i+j}\alpha(\lie{C_i}{C_j},C_0,\dots,\widehat{C_i},\dots,
\widehat{C_j},\dots,C_n) ,
\end{multline*}
satisfies $d_{A\bowtie B}^2=0$.
Now, remember that
\[ \wedge^{n} (A\oplus B)^* = \bigoplus_{k+l=n}
\wedge^k A^*\otimes \wedge^l B^* .\]
It is easy  to see that
\begin{multline*}
d_{A\bowtie B}(\gm (\wedge^k A^*\otimes\wedge^l B^*))
\subset \gm (\wedge^{k+2} A^*\otimes\wedge^{l-1} B^*)
\oplus \gm (\wedge^{k+1} A^*\otimes\wedge^{l} B^*) \\
\oplus \gm (\wedge^{k} A^*\otimes\wedge^{l+1} B^*)
\oplus \gm (\wedge^{k-1} A^*\otimes\wedge^{l+2} B^*) .
\end{multline*}
Moreover, since $A$ and $B$ are Lie subalgebroids of $A\bowtie B$,
 the stronger relation
\[ d_{A\bowtie B} \gm (\wedge^k A^*\otimes\wedge^l B^*)
\subset  \gm  (\wedge^{k+1} A^*\otimes\wedge^{l} B^*)
\oplus \gm (\wedge^{k} A^*\otimes\wedge^{l+1} B^*) \] holds.
Composing $d_{A\bowtie B}$ with the natural projections on each of the direct summands, we get the commutative diagram below:
\[ \xymatrix{
&\gm ( \wedge^k A^*\otimes\wedge^l B^*) \ar[d]^{d_{A\bowtie B}}
\ar[dl]_{\partial_A} \ar[dr]^{(-1)^k\partial_B} & \\
\gm (\wedge^{k+1}A^*\otimes\wedge^l B^*) &
\gm \big((\wedge^{k+1}A^*\otimes\wedge^l B^*) \oplus
(\wedge^k A^*\otimes\wedge^{l+1} B^*) \big)
\ar[l] \ar[r] & \gm (\wedge^k A^*\otimes\wedge^{l+1} B^* ) ,
} \]
where $\partial_A$ and $\partial_B$ are the coboundary operators given by
\eqref{partialA} and \eqref{partialB}.
>From $d_{A\bowtie B}^2=0$, it follows that $\partial_A^2=0$, $\partial_B^2=0$ and $\partial_A\rond\partial_B=\partial_B\rond\partial_A$.

\fboxr{$\Leftarrow$}
Conversely, given the commutative diagram \eqref{eq:double},
one can define an operator
\[ \sections{\wedge^{\bullet}(A\oplus B)^*}\xrightarrow{d_{A\bowtie B}}
\sections{\wedge^{\bullet+1}(A\oplus B)^*} ,\]
whose restriction to $ \gm (\wedge^k A^*\otimes\wedge^l B^* )$ is
$\partial_A+(-1)^k\partial_B$.
Clearly, $d_{A\bowtie B}^2=0$ and $(\sections{\wedge^{\bullet}(A\oplus B)^*}, \ d_{A\bowtie B} )$ is a differential graded algebra.
 Hence it follows that  $A\oplus B$ admits a Lie algebroid structure with associated differential $d_{A\bowtie B}$.
Moreover,
\[ d_{A\bowtie B}\gm (\wedge^k A^*\otimes\wedge^l B^*) 
\subset \gm ( \wedge^{k+1} A^*\otimes\wedge^{l} B^*)
\oplus \gm (\wedge^{k} A^*\otimes\wedge^{l+1} B^*) .\]
Therefore, $\sections{A}$ and $\sections{B}$ are closed under the Lie algebroid bracket on $A\oplus B$. The subbundles $A$ and $B$ are thus Lie subalgebroids of $A\oplus B$. We conclude that the pair $(A,B)$ is a matched pair of Lie algebroids.
\end{proof}

\begin{prop}\label{cohobowtie}
The Lie algebroid cohomology of $A\bowtie B$ (with trivial
coefficients) is isomorphic to the total cohomology of the double
complex \eqref{eq:double}.
\end{prop}

\begin{proof}
This is an immediate consequence of the following fact, which was already pointed out in the proof of Proposition~\ref{doublecomplex}: the restriction of
the cohomology operator
\[ \sections{\wedge^{\bullet}(A\oplus B)^*}\xrightarrow{d_{A\bowtie B}}
\sections{\wedge^{\bullet+1}(A\oplus B)^*} \]
to the subspace $\gm (\wedge^k A^*\otimes\wedge^l B^* )$ of
 $\gm (\wedge^{k+l}(A\oplus B)^* )$ is $\partial_A+(-1)^k\partial_B$.
\end{proof}

\subsection{Canonical complex Lie algebroid associated to a holomorphic Lie algebroid}

The following result is standard, for instance, see \cite{MR2093043, Block2007}.

\begin{lem}\label{lem:Huybrechts}
Let $E$ be a complex vector bundle over a complex manifold $X$.
Then $E$ is a holomorphic vector bundle if, and only if, $E$ is a $T^{0,1} X$-module --- i.e. there exists a \emph{flat} $T^{0,1}X$-connection on $E$:
\[ \sections{T^{0,1} X}\otimes\sections{E}\to\sections{E}:
(X,\varepsilon)\mapsto\nabla_X\varepsilon .\]
\end{lem}

\begin{proof}
\fboxr{$\Rightarrow$}
Let $\shs{E}$ denote the sheaf of holomorphic sections of $E\to X$.
For all $U\subset X$ open and $\sigma\in\shs{E}$(U), set
$\nabla_X\sigma=0,\;\forall X\in\sections{T^{0,1} X|_U}$.
Then $\nabla$ extends to all smooth sections of $E$ by
\[ \nabla_X(f\varepsilon)=X(f)\varepsilon+f\nabla_X\varepsilon \]
since $\sections{E|_U}$ is generated by $\shs{E}(U)$ over $\cinf{U,\CC}$.
One easily sees that $\nabla$ is a flat $T^{0,1} X$-connection.

\fboxr{$\Leftarrow$}
Let $\nabla$ denote the representation of $T^{0,1} X$ on $E$.
And define the sheaf $\shs{E}$ over $X$ by
\[ \shs{E}(U)=\genrel{\sigma\in\sections{E|_U}}
{\nabla_X\sigma=0,\forall X\in\sections{T^{0,1} X|_U}} ,\]
for all $U\subset X$ open.
If $\sigma\in\shs{E}(U)$, then
\[ \nabla_X(f\sigma)=X(f)\sigma+f\nabla_X\sigma=X(f)\sigma ,\]
for all $X\in\sections{T^{0,1} X|_U}$.
Hence $f\sigma\in\shs{E}(U)$ if $f\in\hf{X}(U)$.
Thus $\shs{E}$ is a sheaf of $\hf{X}$-modules.
Since $\nabla$ is flat, given any $e\in E_x$, there exists a neighborhood
$U$ of $x$ and a local section $\sigma\in\shs{E}(U)$ through $e$.
Hence there exists a \emph{holomorphic} vector bundle structure on $E$ with $\shs{E}$
as sheaf of holomorphic sections.
\end{proof}

Now we can state the main result of this section.

\begin{thm}\label{HLA=MP}
Let $A$ be a holomorphic Lie algebroid over a complex manifold $X$.
Then the pair $(T^{0,1} X,A^{1,0})$ is naturally a matched pair of complex Lie algebroids.
Conversely, given a complex manifold $X$ and a matched pair $(T^{0,1} X,B)$,
where $B$ is a complex Lie algebroid over $X$ whose anchor takes its values in
$T^{1,0} X$, there exists a holomorphic Lie algebroid $A$ such that $B\isomorphism A^{1,0}$
as complex Lie algebroids.
\end{thm}

\begin{proof}
\fboxr{$\Rightarrow$}
Let $\anchor$ denote the anchor map of $A$.
Since $A$ is a holomorphic vector bundle, by Lemma~\ref{lem:Huybrechts},
the complex vector bundle $A^{1,0}$ is a $T^{0,1} X$-module.
This gives a representation \[ \sections{T^{0,1} X}\otimes\sections{A^{1,0}}\to\sections{A^{1,0}}:
(X,\eta)\mapsto\nabla_X\eta \] of $T^{0,1} X$ on $A^{1,0}$.
On the other hand, the $\CC$-linear extension of the anchor map
$\anchor: A\to TX$ induces a complex vector bundle map
$\anchor: A^{1, 0}\to T^{1, 0}X$, which is a morphism of complex Lie algebroids.
Similar to the situation of Example~\ref{t01t10}, the map
\begin{equation}
\label{eq:etaX}
 \sections{A^{1,0}}\otimes\sections{T^{0,1} X}\to\sections{T^{0,1} X}:
(\eta,X)\mapsto
%\anchor(\nabla_X\eta)-\lie{X}{\anchor(\eta)}=
\nabla_\eta X:=\pr^{0,1}\lie{\canchor(\eta)}{X} 
\end{equation}
is automatically a representation of $A^{1,0}$ on $T^{0,1} X$.

It remains to prove that the pair $(T^{0,1} X,A^{1,0})$,
with the above two representations, satisfies
the matched pair axioms \eqref{first} to \eqref{third}.

If $\eta$ is a holomorphic section of $A^{1,0}|_U$, then
$\nabla_X\eta=0$ for all $X\in\sections{T^{0,1} X|_U}$
(by definition of the $T^{0,1} X$-module structure of $A^{1,0}$) and,
since $\rho\cmplx(\eta)$ is a holomorphic section of $T^{1,0} X$ over $U$,
we have
\[ \pr^{1,0}\lie{X}{\rho\cmplx(\eta)}=0, \qquad\forall X\in\sections{T^{0,1} X|_U} .\]
Thus
\beq{etiole} \lie{X}{\rho\cmplx(\eta)}=-\pr^{0,1}\lie{\rho\cmplx(\eta)}{X}+\rho\cmplx(\nabla_X\eta) ,\eeq
for all $\eta\in\shs{A}^{1,0} (U)$ and $X\in\sections{T^{0,1} X}$.
>From Lemma~\ref{tensoriality}, it follows that Eq.~\eqref{etiole}
 actually holds for all \emph{smooth} sections $\eta$ of
$A^{1,0}$.
This relation is nothing but axiom~\eqref{first} in the particular case $(T^{0,1} X,A^{1,0})$.

If the sections of $A^{1,0}$ involved in Eq.~\eqref{second} are taken holomorphic,
then \eqref{second} holds because, in that particular case, all its terms vanish.
Therefore, by Lemma~\ref{tensoriality} and the fact that $\shs{A}^{1,0}_\infty$ is generated by $\shs{A}^{1,0}$
as a sheaf of modules over the sheaf $\cinf{X}$ of smooth functions,
Eq.~\eqref{second} is always satisfied.

Finally, it follows from the definition of the $A^{1,0}$-module structure on $T^{0,1} X$ and the Jacobi identity that
\[ \nabla_\eta\lie{X_1}{X_2}=\lie{\nabla_\eta X_1}{X_2}+\lie{X_1}{\nabla_\eta X_2}, \qquad
\forall \eta\in\shs{A}^{1,0}(U),\;\forall X_1,X_2\in\sections{T^{0,1} X} .\]
Hence \eqref{third} follows from Lemma~\ref{tensoriality}.

\fboxr{$\Leftarrow$}
Let $E\to X$ denote the underlying complex vector
 bundle, $\lie{\cdot}{\cdot}_B$
the Lie bracket on $\sections{E}$
and $\anchor:E\to T^{1,0}X$ the anchor map of the Lie
algebroid $B$.
Since $B$ is a $T^{0,1} X$-module, it follows from Lemma~\ref{lem:Huybrechts}
that $E$ is a holomorphic vector bundle --- a smooth
section $\eta \in \sections{E|_U}$ being holomorphic iff
$\nabla_X\eta=0,\;\forall X\in\sections{T^{0,1} X|_U}$.

Moreover, by Eq.~\eqref{second},
if $\eta_1, \ \eta_2\in \sections{E|_U}$
are two holomorphic sections over an
open subset $U\subset X$,
$[\eta_1, \eta_2]$ is also a holomorphic section
of $E$ over $U$; i.e. the sheaf $\shs{E}$ of holomorphic sections of $E$
is a subsheaf of complex Lie subalgebras of the sheaf $\shs{E}_\infty$ of smooth sections.

We define a new Lie algebroid structure $A$ on the vector bundle
 $E$ with the composition \[ E\xrightarrow{\anchor}T^{1,0} X\xrightarrow{2\Re}TX \]
as anchor map and such that the Lie brackets of $A$ and $B$ agree
on the subsheaf $\shs{E}$ of $\shs{E}_\infty$:
\[ \lie{\sigma}{\tau}_A=\lie{\sigma}{\tau}_B
\qquad \forall\sigma,\tau\in\shs{E} (U) .\]
Here $\Re$ means real part.

Eq.~\eqref{first} implies that
\[ \nabla_\eta Y =[\anchor\cmplx(\eta),Y], \qquad\forall\eta\in\shs{E}(U),\;
\forall Y\in\sections{T^{0,1} X|_U} .\] Thus
$\pr^{1,0}\lie{\rho\cmplx(\eta)}{X}=0$. By Example~\ref{t01t10},
$\rho\cmplx(\eta)$ is thus a holomorphic section of $T^{1,0} X|_U$ if
$\eta$ is a holomorphic section of $E|_U$. Hence $\rho : E\to T^{0,1} X$
and $\Re\rond \rho:  E\to TX$ are holomorphic
bundle maps.

Note that
\[
%\begin{multline*}
f\in\hf{X} (U) 
%\quad\Leftrightarrow\quad df\rond J=i\;df
%\quad\Leftrightarrow\quad df\rond\tfrac{1-iJ}{2}=df \\
\quad\Leftrightarrow\quad df(X)=2df(\Re X),\;\forall X\in\sections{T^{1,0}
X} .
%\end{multline*}
\]
Therefore, the Lie bracket on $\shs{A}$ satisfies the Leibniz rule.
 Indeed, for all $\sigma,\tau\in\shs{E}(U)$ and
$f\in\hf{X}(U)$, we have
\[ \lie{\sigma}{f\tau}_A=
\lie{\sigma}{f\tau}_B
= df\big(\rho(\sigma)\big)\tau+f\lie{\sigma}{\tau}_B
= df\big(2\Re\rond\rho(\sigma)\big)\tau+f\lie{\sigma}{\tau}_A .\]

Clearly, $A$ is a holomorphic Lie algebroid over $X$
with the same underlying complex vector bundle $E$
and with $\shs{E}$ as its sheaf of holomorphic sections.
By construction, $A^{1,0}$ and $B$ are isomorphic complex Lie algebroids.
\end{proof}

Thus, given a holomorphic Lie algebroid $A\to X$, we obtain two complex
Lie algebroids: $A\cmplx$ and $T^{0,1}X\bowtie A^{1,0}$. The following
proposition follows easily from the construction of
$T^{0,1}X\bowtie A^{1,0}$.

\begin{prop}
\label{prop:5.9}
Assume $(A\to X,\rho,\lie{\cdot}{\cdot})$ is a holomorphic Lie algebroid.
Then \[ A\cmplx\to T^{0,1}X\bowtie A^{1,0} : a\mapsto\big(\anchor\cmplx(\pr^{0,1}(a)),\pr^{1,0}(a)\big) \]
is a homomorphism of complex Lie algebroids.
\end{prop}

\subsection{Lie algebroid cohomology}

We use the following definition of Lie algebroid cohomology
due to Evens-Lu-Weinstein \cite{MR1726784}.

\begin{defn}
\label{defn:5.11}
Given a \emph{holomorphic} Lie algebroid $A\to X$, let $\Omega^k_A$ be the
sheaf of \emph{holomorphic} sections of $\wedge^k A^*\to X$ ($k=1,2,\dots$) and $\Omega^0_A=\hf{X}$.
We have the following  complex of sheaves over $X$:
\[ \Omega^\com_A: \; \Omega^0_A \xrightarrow{d_A} \Omega^1_A \xrightarrow{d_A}
\cdots \xrightarrow{d_A} \Omega^k_A \xrightarrow{d_A} \Omega^{k+1}_A \xrightarrow{d_A} \cdots \]
where
\begin{align*} \big(d_A\alpha\big)(V_0,\cdots,V_k)=
& \sum_{i=0}^k (-1)^i \rho(V_i) \alpha(V_0,\cdots,\widehat{V_i},\cdots,V_k) \\
& +\sum_{i<j} (-1)^{i+j} \alpha(\lie{V_i}{V_j},V_0,\cdots,\widehat{V_i},\cdots,\widehat{V_j},\cdots,V_k)
\end{align*}
for all $\alpha\in\Omega^k_A(U)$,
$V_0,\cdots,V_k\in\shs{A}(U)$,
and any open subset $U$ of $X$.

The \emph{holomorphic} Lie algebroid cohomology of $A$ (with
trivial coefficients) is defined to be
the cohomology of the complex of sheaves $\Omega^\com_A$:
\[ H^*(A,\CC) := H^*(X,\Omega^\com_A) .\]
\end{defn}

The following result gives us an important way of
computing holomorphic Lie algebroid cohomology using
smooth cohomology, i.e. Lie algebroid cohomology
of complex Lie algebroids. In a certain
sense, this is a generalization of Dolbeault's theorem
to Lie algebroids.

\begin{thm}
\label{thm:5.8}
Let $A\to X$ be a \emph{holomorphic} Lie algebroid. Then
\[ H^*(A,\CC) \isomorphism H^*(T^{0,1}X\bowtie A^{1,0},\CC) ,\]
where the right hand side stands for the \emph{complex} Lie algebroid
cohomology of $T^{0,1}X\bowtie A^{1,0}$ (see Proposition~\ref{cohobowtie}), 
which can be interpreted as a generalization of the Dolbeault cohomology.
\end{thm}

\begin{proof}
By $\OA{k}{l}$, we denote the
sheaf of sections of the complex vector bundle
$(T^{0, k}X)^* \otimes \wedge^l A^{1, 0} \to X$.
By the holomorphic Poincar\'e lemma, we have the following
 resolution of complex of sheaves: 
\begin{equation}\label{eq:A}
\xymatrix{
\cdots&\cdots&\cdots&\cdots\\
0 \ar[r] &\Omega_A^2 \ar[u]^{d_A} \ar[r]^\delbar &
\OA{0}{2}\ar[u]^{d_A^{1, 0}}\ar[r]^\delbar
&\OA{1}{2} \ar[u]^{d_A^{1, 0}}\ar[r]^\delbar
&\cdots\\
0 \ar[r] &\Omega_A^1 \ar[u]^{d_A} \ar[r]^\delbar &
\OA{0}{1}\ar[u]^{d_A^{1, 0}}\ar[r]^\delbar
&\OA{1}{1} \ar[u]^{d_A^{1, 0}}\ar[r]^\delbar
&\cdots\\
0 \ar[r] &\Omega_A^0\ar[u]^{d_A}\ar[r]^\delbar &
\OA{0}{0}\ar[u]^{d_A^{1, 0}}\ar[r]^\delbar
& \OA{1}{0}  \ar[u]^{d_A^{1, 0}}\ar[r]^\delbar &\cdots
}
\end{equation}
The conclusion thus follows immediately from Theorem~\ref{HLA=MP}.
\end{proof}

\begin{ex}
As in Example~\ref{t01t10}, consider a complex manifold $X$.
Let $A=TX$ be its holomorphic tangent bundle considered
as a holomorphic Lie algebroid. In this particular case, $H^*(A,\CC)$
is the holomorphic de~Rham cohomology, while
$H^*(T^{0,1}X\bowtie A^{1,0},\CC)$ is the 
$\CC$-valued smooth de~Rham cohomology  since $T^{0,1}X\bowtie A^{1,0}
=T^{0,1} X\bowtie T^{1,0} X\isomorphism T_\CC X$ as complex Lie algebroids.
It is well known that they are isomorphic.
\end{ex}

Indeed the complex Lie algebroid $T^{0,1} X\bowtie A^{1, 0}$ is
an elliptic Lie algebroid in the sense of Block. 
Recall that a complex Lie algebroid $(A\to X, [\cdot, \cdot], \rho)$
is said to be {\em elliptic} \cite{Block2007} if
the map $\Re\rond \rho: A\to TX$ is surjective.
For an elliptic  Lie algebroid, the Lie algebroid cohomology 
cochain complex is an elliptic complex \cite{Block2007}.
Hence the cohomology groups are finite dimensional
if the base manifold is compact. This can be easily seen directly from our
definition of holomorphic Lie algebroid cohomology in terms of a complex of sheaves.

Denote by $H^*(A_R,\CC)$ the Lie algebroid cohomology of the
underlying real Lie algebroid $A_R$ with trivial complex coefficients. 
The following result is an immediate consequence of
Theorem~\ref{thm:5.8} and
Proposition~\ref{prop:5.9}.

\begin{prop}\label{prop:hs}
Let $A$ be a holomorphic Lie algebroid with
underlying real Lie algebroid $A_R$. Then there
is a natural morphism
\[ H^*(A,\CC)\to H^*(A_R,\CC) .\]
\end{prop}

\begin{rmk}
In \cite{MR2285039}, Weinstein asked the question how
to integrate a complex Lie algebroid. For 
the complex Lie algebroid $T^{0,1}X\bowtie A^{1,0}$
arising from a matched pair, our discussion above
suggests that the holomorphic Lie groupoid integrating the
holomorphic Lie algebroid $A$ might be a candidate. It will
be interesting to explore further if Theorem \ref{HLA=MP}
would have any applications in solving Weinstein's integration
problem.
\end{rmk}

\subsection{Cohomology with general coefficients}

\begin{defn}
Given a holomorphic Lie algebroid $A\to X$, an $A$-module is
a holomorphic vector bundle $E\to X$
together with a morphism of sheaves (of $\CC$-modules)
\[ \shs{A}\otimes\shs{E}\to\shs{E}:V\otimes s\mapsto\nabla_V s \]
such that, for any open subset $U\subset X$, the relations
\begin{gather*}
\nabla_{fV} s=f\nabla_V s \\
\nabla_V (fs)=\big(\rho(V)f\big) s +f \nabla_V s \\
\nabla_V\nabla_W s-\nabla_W\nabla_V s=\nabla_{\lie{V}{W}} s
\end{gather*}
are satisfied $\forall f\in\hf{X}(U)$, $\forall V,W\in\shs{A}(U)$
and $\forall s\in\shs{E}(U)$.
\end{defn}

\begin{lem}
Let $A\to X$ be a \emph{holomorphic} Lie algebroid and $E\to X$ a complex
vector bundle. Then $E\to X$ is an $A$-module if, and only if,
$E\to X$ is a module over the complex Lie algebroid
$T^{0,1}X\bowtie A^{1,0}$.
\end{lem}

\begin{proof}
\fboxr{$\Rightarrow$}
Firstly, note that, since $E\to X$ is a holomorphic vector bundle with
sheaf of holomorphic sections $\shs{E}$, by Lemma~\ref{lem:Huybrechts},
$E$ is a $T^{0,1} X$-module whose representation map
\[ \sections{T^{0,1} X}\otimes\sections{E}\to\sections{E}:(X,s)\mapsto\nabla^{0,1}_X s \]
is entirely characterized by the condition
\[ \shs{E}(U)=\genrel{\sigma\in\sections{E|_U}}{\nabla^{0,1}_X\sigma=0,\;\forall X\in\sections{T^{0,1} X|_U}} ,\]
for all open subsets $U\in X$.

On the other hand, $E$ is a module over the holomorphic Lie algebroid $(A\to X,\lie{\cdot}{\cdot},\rho)$
with representation map
\[ \shs{A}\otimes\shs{E}\to\shs{E}:(a,s)\mapsto\nabla_a s ,\]
where $\shs{A}$ and $\shs{E}$ are the $\hf{X}$-sheaves of holomorphic sections of $A$ and $E$ respectively.
Since $\shs{A}_\infty$ and $\shs{E}_\infty$ are generated by $\shs{A}$ and 
$\shs{E}$ as sheaves of $C_X^\infty$-modules,
one can define a representation
\[ \sections{A^{1,0}}\otimes\sections{E}\to\sections{E}:(a,s)\mapsto\nabla^{1,0}_a s \]
of the complex Lie algebroid $(A^{1,0},\lie{\cdot}{\cdot}^{1,0},\rho^{1,0})$ on $E$,
which is entirely determined by the requirement that, for all subsets $U$ of $X$, one has
\[ \nabla^{1,0}_{(1-ij)\alpha}\sigma=\nabla_\alpha \sigma \]
for all $\alpha\in\shs{A}(U)$ and $\sigma\in\shs{E}(U)$.

The equation \[ \nabla_{(X,a)}s=\nabla^{0,1}_X s+\nabla^{1,0}_a s \] defines
a connection of the complex Lie algebroid $T^{0,1} X\bowtie A^{1,0}$
(associated to the holomorphic Lie algebroid $A$ as in 
Theorem~\ref{HLA=MP}) on $E$.
To check that this covariant derivative is flat, it suffices to prove that
\[ \nabla_{(X,0)}\nabla_{(0,a)}s-\nabla_{(0,a)}\nabla_{(X,0)}s=\nabla_{\lie{(X,0)}{(0,a)}}s ,\]
for all $X\in\sections{T^{0,1} X}$, $a\in\sections{A^{1,0}}$ and $s\in\sections{E}$.
However, the curvature being a tensor, it actually suffices to check that,
for any open subset $U$ of $X$, one has
\beq{star27} \nabla_{(X,0)}\nabla_{(0,(1-ij)\alpha)}\sigma-\nabla_{(0,(1-ij)\alpha)}\nabla_{(X,0)}\sigma
=\nabla_{\lie{(X,0)}{(0,(1-ij)\alpha)}}\sigma \eeq
for all $X\in\sections{T^{0,1} X|_U}$, $\alpha\in\shs{A}(U)$ and $\sigma\in\shs{E}(U)$.
By definition of the Lie algebroid structure of $T^{0,1} X\bowtie A^{1,0}$,
\[ \lie{(X,0)}{(0,(1-ij)\alpha)}=(-\nabla_\alpha X,0) \]
since $\alpha$ is holomorphic.
Hence, Eq.~\eqref{star27} becomes
\[ \nabla^{0,1}_X\nabla_\alpha \sigma-\nabla_\alpha\nabla^{0,1}_X\sigma=\nabla^{0,1}_{-\nabla_\alpha X}\sigma ,\]
which is obviously true since each term in the equation above  vanishes
--- indeed, the second argument of each $\nabla^{0,1}$ is a holomorphic section of $E$.

\fboxr{$\Leftarrow$}
Let $\nabla$ denote the representation of $T^{0,1}X\bowtie A^{1,0}$ on $E$.

Since $T^{0,1}X$ is a Lie subalgebroid of $T^{0,1}X\bowtie A^{1,0}$, $E$ is a $T^{0,1}X$-module
and thus, by Lemma~\ref{lem:Huybrechts}, $E\to X$ is a holomorphic vector bundle whose sheaf
of holomorphic sections $\shs{E}$ is characterized by
\beq{star1} \shs{E}(U)=\genrel{\sigma\in\sections{E|_U}}{\nabla_{(X,0)}\sigma=0,
\;\forall X\in\sections{T^{0,1} X|_U}} \eeq

Moreover, the curvature of $\nabla$ being flat, we have
\beq{star0} \nabla_{(X,0)}\nabla_{(0,(1-ij)\alpha)}\sigma-\nabla_{(0,(1-ij)\alpha)}\nabla_{(X,0)}\sigma
=\nabla_{\lie{(X,0)}{(0,(1-ij)\alpha)}}\sigma ,\eeq
for any open subset $U$ of $X$ and all $X\in\sections{T^{0,1} X|_U}$, $\alpha\in\shs{A}(U)$ and $\sigma\in\shs{E}(U)$.
Note that, since $\alpha$ is holomorphic,
\beq{star2} \lie{(X,0)}{(0,(1-ij)\alpha)} \in \sections{T^{0,1} X|_U} .\eeq

Making use of Eqs.~\eqref{star1}~and~\eqref{star2}, 
Eq.~\eqref{star0} becomes
\[ \nabla_{(X,0)} \big(\nabla_{(0,(1-ij)\alpha)}\sigma\big) = 0 ,\]
for any open subset $U$ of $X$ and all $X\in\sections{T^{0,1} X|_U}$, $\alpha\in\shs{A}(U)$ and $\sigma\in\shs{E}(U)$.
In other words, the map
\[ \shs{A}(U)\times\shs{E}(U)\to\shs{E}(U):(\alpha,\sigma)\mapsto\nabla_{(0,(1-ij)\alpha)}\sigma \]
does indeed take its values in $\shs{E}$.
Therefore this restriction of $\nabla$ endows the holomorphic vector bundle $E\to X$ with a structure of module
over the holomorphic Lie algebroid $A$.
\end{proof}

\begin{defn}
Given a \emph{holomorphic} Lie algebroid $A\to X$
and an $A$-module $(E\to X,\nabla)$, we have the complex of sheaves over $X$
\[ \Omega^\com_A\otimes\shs{E}: \; \Omega^0_A\otimes\shs{E}
\xrightarrow{d_A^\nabla} \Omega^1_A\otimes\shs{E} \xrightarrow{d_A^\nabla}
\cdots \xrightarrow{d_A^\nabla} \Omega^k_A\otimes\shs{E} \xrightarrow{d_A^\nabla}
\Omega^{k+1}_A\otimes\shs{E} \xrightarrow{d_A^\nabla} \cdots \]
where
\begin{align*} \big(d_A^\nabla\alpha\big)(V_0,\cdots,V_k)=
& \sum_{i=0}^k (-1)^i \nabla_{V_i} \big(\alpha(V_0,\cdots,\widehat{V_i},\cdots,V_k)\big) \\
& +\sum_{i<j} (-1)^{i+j} \alpha(\lie{V_i}{V_j},V_0,\cdots,\widehat{V_i},\cdots,\widehat{V_j},\cdots,V_k)
\end{align*}
for any open subset $U$ of $X$ and all
$\alpha\in\big(\Omega^k_A\otimes\shs{E}\big)(U)$
and $V_0,\cdots,V_k\in\shs{A}(U)$.

The \emph{holomorphic} Lie algebroid cohomology of $A$ (with
coefficients in the $A$-module $E$) is defined to be
the cohomology of the complex of sheaves $\Omega^\com_A\otimes\shs{E}$:
\[ H^*(A,E) := H^*(X,\Omega^\com_A\otimes\shs{E}) .\]
\end{defn}

\begin{lem}
Let $A\to X$ be a \emph{holomorphic} Lie algebroid and $E\to X$
an $A$-module.
If   $X$ is  compact, then $H^k(A, E)$ is finite 
dimensional for all $k$.
\end{lem}

The following theorem can be proved in a similar fashion as in
Theorem~\ref{thm:5.8}.

\begin{thm}\label{thm:5.9}
Let $A\to X$ be a \emph{holomorphic} Lie algebroid and $E\to X$
an $A$-module. Then
\[ H^*(A, E) \isomorphism H^*(T^{0,1}X\bowtie A^{1,0}, E) .\]
\end{thm}

Given a holomorphic Lie algebroid $A\to X$  and an $A$-module 
$E\to X$, it is simple to see that
$E\to X$ naturally becomes an $A_R$-module.
The following is a straightforward generalization of Proposition~\ref{prop:hs}.

\begin{prop}
Let $A$ be  a holomorphic Lie algebroid with
underlying real Lie algebroid $A_R$, and
$E\to X$ an $A$-module. Then there
is a natural homomorphism \[ H^*(A,E)\to H^*(A_R,E) .\]
\end{prop}

\subsection{Application to holomorphic Poisson manifolds}

Now consider the cotangent bundle Lie algebroid  $(T^*X)_\pi$ associated
to a holomorphic Poisson structure $(X, \pi)$.
Since $(T^*X)_\pi$ is a holomorphic Lie algebroid, according
to Theorem~\ref{HLA=MP},  $(T^{0,1}X,  (T^{1, 0} X)^*_\pi )$
is a matched pair.
On the other hand, according to Theorem~\ref{PNGC},
 to any holomorphic Poisson manifold
corresponds a natural generalized complex structure.
The following theorem indicates the relation between
this generalized complex structure and the
complex Lie algebroid  $T^{0,1}X\bowtie (T^{1, 0} X)^*_\pi$.

\begin{thm}\label{YAO}
If $(X, \pi)$ is a holomorphic Poisson manifold, then
the  complex Lie algebroid $T^{0, 1}X\bowtie (T^{1, 0} X)^*_\pi$ is
isomorphic to the Dirac structure $L_{4\pi}$,  the $-i$-eigenbundle
 of the generalized complex structure
\[ \JJ_{4\pi}= \begin{pmatrix} J & 4\piim\diese \\ 0 & -J^* \end{pmatrix} \]
as in Theorem~\ref{PNGC}.
\end{thm}

We need a few lemmas.

\begin{lem}
\begin{equation}
L_{4\pi}=\{(X^{0,1}+\pi\diese\xi^{1, 0}, \xi^{1, 0})
|\xi^{1, 0}\in \Omega^{1, 0}(X), \; X^{0,1}\in \XX^{0,1} (X)\}
\end{equation}
\end{lem}

\begin{proof}
It is clear that \[ \JJ_{4\pi}(X^{0,1}, 0)=(JX^{0,1}, 0)=-i
(X^{0,1}, 0) .\] On the other hand, since
$\piim\diese=\tfrac{1}{2i}(\pi\diese-\overline{\pi}\diese)$, it follows
that
\[ \piim\diese\xi^{1, 0}=\tfrac{1}{2i}(\pi\diese\xi^{1, 0})
=-\tfrac{i}{2}\pi\diese \xi^{1, 0} \] Since
$\pire\diese=J\rond\piim\diese$, we have
\[ J\rond \pi\diese=J\rond (\pire\diese+ i\piim\diese )=-\piim\diese +i\pire\diese
=i (\pire\diese+ i\piim\diese )=i\pi\diese .\]
It thus follows that
\[ \JJ_{4\pi} (\pi\diese\xi^{1, 0}, \xi^{1, 0})=(J \pi\diese\xi^{1, 0}+4\piim\diese\xi^{1,0}, -J^*\xi^{1, 0})
=(-i \pi\diese \xi^{1, 0}, -i\xi^{1, 0})=-i (\pi\diese \xi^{1, 0}, \xi^{1, 0}) .\]
Hence $(X^{0,1}+\pi\diese\xi^{1, 0}, \xi^{1, 0})$ is an  eigenvector
of $\JJ_{4\pi}$ with eigenvalue $-i$. The conclusion thus follows
from dimension counting.
\end{proof}

By abuse of notations, $\nabla$ denotes both the 
$T^{0,1}X$-representation on $(T^{1, 0} X)^*_\pi$
and the $(T^{1, 0} X)^*_\pi$-representation on
$T^{0,1}X$.

\begin{lem}
For any $X^{0, 1}\in \XX^{0, 1}(X)$ and
$\xi^{1, 0}\in \Omega^{1, 0}(X)$, we have
\begin{equation}\label{35}
\nabla_{X^{0, 1}}\xi^{1, 0}=L_{X^{0, 1}}\xi^{1, 0}
.\end{equation}
\end{lem}

\begin{proof}
For any $Y^{1, 0}\in \XX^{1, 0}(X)$, we have
\begin{align*}
\pairing{\nabla_{X^{0,1}}\xi^{1,0}}{Y^{1,0}}
=& X^{0,1}\pairing{\xi^{1,0}}{Y^{1,0}} - \pairing{\xi^{1,0}}{\nabla_{X^{0,1}}Y^{1,0}} \\
=& X^{0,1}\pairing{\xi^{1,0}}{Y^{1,0}} - \pairing{\xi^{1,0}}{\pr^{1,0}\lie{X^{0,1}}{Y^{1,0}}} \\
=& X^{0,1}\pairing{\xi^{1,0}}{Y^{1,0}} - \pairing{\xi^{1,0}}{\lie{X^{0,1}}{Y^{1,0}}} \\
=& \pairing{L_{X^{0,1}}\xi^{1,0}}{Y^{1,0}}
.\end{align*}
Hence $\nabla_{X^{0,1}}\xi^{1,0}=L_{X^{0,1}}\xi^{1,0}$.
\end{proof}

\begin{lem}
For any $X^{0, 1}\in \XX^{0, 1}(X)$ and
$\xi^{1, 0}\in \Omega^{1, 0}(X)$, we have
\begin{equation}\label{36}
\pi\diese\nabla_{X^{0,1}}\xi^{1,0}=\pr^{1, 0}[\pi\diese\xi^{1,0},X^{0,1}]
.\end{equation}
\end{lem}

\begin{proof}
Note that if $Y^{1,0}\in\XX^{1,0}(X)$ is a holomorphic vector field,
$\pr^{1,0}[X^{0,1},Y^{1,0}]=0$. Hence it follows that
$L_{X^{0,1}}\pi\in\sections{T^{0,1}X\wedge T^{1,0}X}$, for
$\pi\in\vf^{2,0}(X)$ is a holomorphic bivector field.
Therefore
$(L_{X^{0,1}}\pi)\diese\xi^{1,0}\in\XX^{0,1}(X)$. That is,
$\pr^{1, 0}(L_{X^{0,1}}\pi)\diese\xi^{1, 0}=0$.
Now since
\[ \pi\diese(L_{X^{0,1}}\xi^{1,0})-L_{X^{0,1}}(\pi\diese\xi^{1,0})
=(L_{X^{0,1}}\pi)\diese\xi^{1,0} ,\]
applying  $\pr^{1,0}$ to both sides, we obtain
\[ \pi\diese(L_{X^{0,1}}\xi^{1,0})=\pr^{1,0}[X^{0,1},\pi\diese\xi^{1,0}] \]
and, using Eq.~\eqref{35},
\[ \pi\diese(\nabla_{X^{0,1}}\xi^{1,0})=\pr^{1,0}\lie{X^{0,1}}{\pi\diese\xi^{1,0}} . \qedhere\]
\end{proof}

\begin{proof}[Proof of Theorem~\ref{YAO}]
First, recall that the Lie bracket on $\sections{L_{4\pi}}$ is the restriction of the 
Courant bracket \[ \courant{X+\xi}{Y+\eta}=\lie{X}{Y}+\derlie{X}\eta-\derlie{Y}{\xi}
+\thalf d(\xi Y-\eta X) \] of $\sections{TX\oplus T^*X}$ \cite{Gualtieri2007a,MR2276462}.

Consider the map
\[ \phi: T^{0, 1}X\bowtie (T^{1, 0} X)^*_\pi\to L_{4\pi} :
(X^{0, 1}, \xi^{1, 0} )\mapsto (X^{0, 1}+\pi\diese\xi^{1, 0}, \xi^{1, 0}) , \]
which is an isomorphism of vector bundles.
It is clear that $\phi$ interchanges the anchor maps.
One sees immediately that
\begin{align*}
[\phi(X^{0,1}),\phi(Y^{0,1})] &= \phi[X^{0,1},Y^{0,1}] &&\forall X^{1,0},Y^{1,0}\in\vf^{1,0}(X) \\
[\phi(\xi^{1,0}),\phi(\eta^{1,0})] &= \phi[\xi^{1,0},\eta^{1,0}] &&
\forall \xi^{1,0},\eta^{1,0}\in\Omega^{1,0}(X).
\end{align*}

Now for any $X^{0, 1}\in \XX^{0, 1}(X)$ and
$\xi^{1, 0}\in \Omega^{1, 0}(X)$, 
\begin{align*}
&[\phi (X^{0, 1}), \phi (\xi^{1, 0})] && \\
=& \courant{X^{0,1}}{\pi\diese\xi^{1,0}+\xi^{1,0}}
&&\text{(by definition of } \phi) \\
=& [X^{0,1},\pi\diese\xi^{1,0}]+L_{X^{0,1}}\xi^{1,0}
&&\text{(by definition of }\courant{\cdot}{\cdot}) \\
=& [X^{0,1},\pi\diese\xi^{1,0}]+\nabla_{X^{0,1}}\xi^{1,0}
&&\text{(by Eq.~\eqref{35})} \\
=& \pr^{1,0}[X^{0,1},\pi\diese\xi^{1,0}]+\pr^{0,1}
[X^{0,1},\pi\diese\xi^{1,0}]+\nabla_{X^{0,1}}\xi^{1,0}
&&\text{(since }T\cmplx X=T^{1,0} X\oplus T^{0,1} X\text{)} \\
=& \pr^{1,0}[X^{0,1},\pi\diese\xi^{1,0}]-\nabla_{\xi^{1,0}}X^{0,1}
+\nabla_{X^{0,1}}\xi^{1,0} &&\text{(by Eq.~\eqref{eq:etaX})}
\end{align*}
On the other hand,
\begin{multline*}
\phi[X^{0,1},\xi^{1,0}] = \phi(-\nabla_{\xi^{1,0}}X^{0,1}
+\nabla_{X^{0,1}}\xi^{1,0}) %\\
=-\nabla_{\xi^{1,0}}X^{0,1}+\pi\diese\nabla_{X^{0,1}}\xi^{1,0}
+\nabla_{X^{0,1}}\xi^{1,0} .\end{multline*} From Eq.~\eqref{36},
it thus follows that
\[ [\phi(X^{0,1}),\phi(\xi^{1,0})]=\phi[X^{0,1},\xi^{1,0}] .\]
Hence $\phi$ is indeed a Lie algebroid isomorphism.
\end{proof}

\subsection{Holomorphic Poisson cohomology}

Consider the matched pair $(T^{0,1}X, (T^{1,0}X)^*_\pi)$
associated to a holomorphic Poisson manifold $(X, \pi)$. With 
Proposition~\ref{doublecomplex} in mind, we set $A =
T^{0,1}X $ and $ B=(T^{1,0}X)_\pi^* $. Then
\[ \Gamma (\wedge^k A^* \otimes  \wedge^l B^*) \simeq
\Omega^{0,k} (X)  \otimes_{{\mathcal C}^{ \infty } (X, \CC) }
{\mathfrak X}^{l,0} (X) \simeq \Omega^{0,k} ( X, T^{l,0}X ) \]
and the commutative diagram of Proposition~\ref{doublecomplex} becomes
\beq{eq:doublexxx} \xymatrix{ \Omega^{0,k} ( X, T^{l,0}X )
\ar[r]^{\partial_A} \ar[d]_{\partial_B} &
\Omega^{0,k+1} ( X, T^{l,0}X ) \ar[d]^{\partial_B} \\
\Omega^{0,k} ( X, T^{l+1,0}X ) \ar[r]_{\partial_A} & \Omega^{0,k+1} (
X, T^{l+1,0}X ). }
 \eeq

The next proposition describes the coboundary operators
$\partial_A$ and $\partial_B$ in this context.

\begin{prop}
Let $(T^{0,1}X,(T^{1,0}X)^*_\pi)$ be the matched pair associated 
to a holomorphic Poisson manifold $(X,\pi)$. Then
\begin{enumerate}
\item $\partial_A:\Omega^{0,k}(X,T^{l,0}X)\to\Omega^{0,k+1}(X,T^{l,0}X)$ 
is the $\delbar$-operator associated to the holomorphic vector bundle $T^{l,0}X$; 
\item $\partial_B:\Omega^{0,k}(X,T^{l,0}X)\to\Omega^{0,k}(X,T^{l+1,0}X)$ 
is the operator $d_\pi$ defined by the relation
\beq{eq:parB}
\big(d_\pi\alpha\big)(Y_1,\dots,Y_k)
=[\pi,\alpha(Y_1,\cdots,Y_k)]+(-1)^k\alpha\per[\pi,Y_1\wedge\dots\wedge Y_k]
,\eeq
where $Y_1,\dots,Y_k$ are arbitrary elements of $\XX^{0,1}(X)$ 
and $\alpha\per[\pi,Y_1\wedge\dots\wedge Y_k]$ denotes
the element in $\XX^{l+1,0}(X)$ obtained by contracting $\alpha$ with
the $(k+1)$-vector field $[\pi,Y_1\wedge\dots\wedge Y_k]$.

Alternatively, if $\omega\in\Omega^{0,k}(X)$ and $P\in\mathfrak{X}^{l,0}(X)$, then 
\begin{equation}\label{eq:parBB}
d_\pi(\omega\otimes P)=\omega\otimes[\pi,P]
+\sum_{i=1}^n(\imath_{\pi\diese(e^i)}\diff\omega)\otimes(e_i\wedge P),
\end{equation}
where $(e_1,\dots,e_n)$ is a basis of $T^{1,0}_x X$ 
and $(e^{1},\dots,e^{n})$ is the dual basis of $(T^{1,0}_x X)^*$.
\end{enumerate}
\end{prop}

\begin{proof}
\fboxr{a} This is straightforward. 
\fboxr{b} 
Since Eq.~\eqref{eq:parB} follows easily from Eq.~\eqref{eq:parBB}, we will only prove Eq.~\eqref{eq:parBB}.

For any $A_1,\cdots,A_k\in\XX^{0,1}(X)$ and
$B_0,\dots,B_l\in\Omega^{1,0}(X)$, according to Eq.~\eqref{partialB}, we have
\begin{multline*}
\big(\partial_B(\omega\otimes P)\big)(A_1,\dots,A_k,B_0,\dots,B_l) 
=T\cdot\omega(A_1,\dots,A_k) \\ 
+ \sum_{i=0}^l (-1)^i S_i\cdot P(B_0,\dots,\widehat{B_i},\dots,B_l) 
\end{multline*}
Here
\begin{multline*} 
T=\sum_{j=0}^l (-1)^j \pi\diese(B_j) P(B_0,\dots,\widehat{B_j},\dots,B_l) \\
+ \sum_{j_1,j_2=0}^l (-1)^{j_1+j_2} P([B_{j_1},B_{j_2}],B_0,\cdots,\widehat{B_{j_1}},\dots, \widehat{B_{j_2}},\dots,B_l) 
\end{multline*}
and, for $i=1,\dots,l$, 
\[ S_i=\pi\diese(B_i)(\omega(A_1,\dots,A_k))
-\sum_{j=1}^k \omega(A_1,\dots,\nabla_{B_i}A_j,\dots,A_k) .\]
It is clear  that
\begin{equation}\label{eq:lieder1}
T=[\pi,P](B_0,\dots,B_l) 
.\end{equation}
According to Eq.~\eqref{eq:etaX}, we have $\nabla_{B_i}A_j=\pr^{0,1}[\pi\diese B_i,A_j]$. 
Since $\omega$ is a $(0,k)$-form, it follows that
\[ ([\pi\diese B_i,A_j]-\pr^{0,1}[\pi\diese B_i,A_j])\per\omega=0 .\]
As a consequence, we have
\[ \omega(A_1,\dots,\nabla_{B_i}A_j,\dots,A_k)
=\omega(A_1,\dots,[\pi\diese B_i,A_j],\dots,A_k) .\]
Therefore
\begin{align*}
S_i=& \pi\diese(B_i)\big(\omega(A_1,\dots,A_k)\big)
-\sum_{j=1}^k\omega(A_1,\dots,[\pi\diese B_i,A_j],\dots,A_k) \\
=& \big(L_{\pi\diese(B_i)}\omega\big)(A_1,\dots,A_k) \\
=& \big(\imath_{\pi\diese(B_i)}\diff\omega\big) (A_1,\dots,A_k) 
,\end{align*}
where the last equality uses the relation $\imath_{\pi\diese B_i}\omega=0$. 
Hence it follows that
\begin{align*}
& \sum_{i=0}^l (-1)^i S_i \cdot P(B_0,\dots,\widehat{B_i},\dots,B_l) \\
=& \sum_{i=0}^l (-1)^i \big(\imath_{\pi\diese(B_i)}\diff\omega\big)(A_1,\dots,A_k) \; P(B_0,\dots,\widehat{B_i},\dots, B_l) \\
=& \sum_{i=0}^l (-1)^i \sum_{j=1}^n B_i(e_j) \;
\big(\imath_{\pi\diese(e^j)}\diff\omega\big)(A_1,\dots,A_k) \;
P(B_0,\dots,\widehat{B_i},\dots,B_l) \\
=& \sum_{j=1}^n \big(\imath_{\pi\diese(e^j)}\diff\omega\big)(A_1,\dots,A_k)
\; \big(e_j\wedge P\big)(B_0,\dots,B_l)
.\end{align*}
This concludes the proof of the proposition.
\end{proof}

As an immediate consequence, we have the following

\begin{cor}
Let $(X,\pi)$ be a holomorphic Poisson manifold. The
following cohomologies are all isomorphic:
\begin{enumerate}
\item the  holomorphic Poisson cohomology of $(X,\pi)$;
\item the complex Lie algebroid cohomology of $L_{4\pi}$;
\item the total cohomology of the double complex
\beq{eq:DeRham}
\xymatrix{ \cdots & \cdots & \cdots & \\
\Omega^{0,0}(X,T^{2,0}X) \ar[u]^{\dpi}\ar[r]^{\delbar} & \Omega^{0,1}(X,T^{2,0}X) \ar[u]^{\dpi}\ar[r]^{\delbar} &
\Omega^{0,2}(X,T^{2,0}X) \ar[u]^{\dpi}\ar[r]^{\delbar} & \cdots \\
\Omega^{0,0}(X,T^{1,0}X) \ar[u]^{\dpi}\ar[r]^{\delbar} & \Omega^{0,1}(X,T^{1,0}X) \ar[u]^{\dpi}\ar[r]^{\delbar} &
\Omega^{0,2}(X,T^{1,0}X) \ar[u]^{\dpi}\ar[r]^{\delbar} & \cdots \\
\Omega^{0,0}(X,T^{0,0}X) \ar[u]^{\dpi}\ar[r]^{\delbar} & \Omega^{0,1}(X,T^{0,0}X) \ar[u]^{\dpi}\ar[r]^{\delbar} &
\Omega^{0,2}(X,T^{0,0}X) \ar[u]^{\dpi}\ar[r]^{\delbar} & \cdots
} \eeq
\end{enumerate}
Here $d_\pi$ is the differential operator defined by Eq.~\eqref{eq:parB} or Eq.~\eqref{eq:parBB}.
\end{cor}

\begin{rmk}
When $X$ is a Stein manifold (for instance $X=\CC^n$),
one easily sees that our Poisson cohomology  groups are isomorphic
to the ones defined by Lichn\'erowicz's
 cochain complex of holomorphic mutivector fields,
as in the smooth case (see, for instance,  \cite{MR2228339}). 

We also note that our holomorphic Poisson cohomology 
groups are always finite dimensional if the manifold is compact.
\end{rmk}

%%%%%%%%%%%%%%%%%%%%%%%%%%%%%%%%%%%%%%%%%%%%%%%%%

\def\cprime{$'$}
\providecommand{\bysame}{\leavevmode\hbox to3em{\hrulefill}\thinspace}
\providecommand{\href}[2]{#2}


\begin{thebibliography}{10}

\bibitem{Block2007}
J.~Block, \emph{Duality and equivalence of module categories in noncommutative
  geometry {I}},  (2007), \mbox{\arxiv{math/0509284}}.

\bibitem{Bondal1993}
A.~Bondal, \emph{Non-commutative deformations and {P}oisson brackets on
  projective spaces}, Preprint no. 67, Max Planck Institut Bonn, 1993.

\bibitem{MR2180064}
M.~N. Boyom, \emph{K{V}-cohomology of {K}oszul-{V}inberg algebroids and
  {P}oisson manifolds}, Internat. J. Math. \textbf{16} (2005), no.~9,
  1033--1061.

\bibitem{MR2263715}
K.~A. Brown, K.~R. Goodearl, and M.~Yakimov, \emph{Poisson structures on affine
  spaces and flag varieties. {I}. {M}atrix affine {P}oisson space}, Adv. Math.
  \textbf{206} (2006), no.~2, 567--629.

\bibitem{MR1665693}
J.~L. Brylinski and G.~Zuckerman, \emph{The outer derivation of a complex
  {P}oisson manifold}, J. Reine Angew. Math. \textbf{506} (1999), 181--189.

\bibitem{MR1718638}
S.~Chemla, \emph{A duality property for complex {L}ie algebroids}, Math. Z.
  \textbf{232} (1999), no.~2, 367--388.

\bibitem{MR996653}
A.~Coste, P.~Dazord, and A.~Weinstein, \emph{Groupo\"\i des symplectiques},
  Publications du D\'epartement de Math\'ematiques. Nouvelle S\'erie. A, Vol.\
  2, Publ. D\'ep. Math. Nouvelle S\'er. A, vol.~87, Univ. Claude-Bernard, 1987.

\bibitem{Crainic2007}
M.~Crainic, \emph{Generalized complex structures and {L}ie brackets},  (2007),
  \mbox{\arxiv{math/0412097}}.

\bibitem{Damianou2007}
P.~A. Damianou and R.~L. Fernandes, \emph{Integrable hierarchies and the
  modular class},  (2007), \mbox{\arxiv{math/0607784}}.

\bibitem{MR725930}
V.~G. Drinfel{\cprime}d, \emph{Constant quasiclassical solutions of the
  {Y}ang-{B}axter quantum equation}, Dokl. Akad. Nauk SSSR \textbf{273} (1983),
  no.~3, 531--535.

\bibitem{MR688240}
\bysame, \emph{Hamiltonian structures on {L}ie groups, {L}ie bialgebras and the
  geometric meaning of classical {Y}ang-{B}axter equations}, Dokl. Akad. Nauk
  SSSR \textbf{268} (1983), no.~2, 285--287.

\bibitem{MR1862022}
S.~Evens and J.-H. Lu, \emph{On the variety of {L}agrangian subalgebras. {I}},
  Ann. Sci. \'Ecole Norm. Sup. (4) \textbf{34} (2001), no.~5, 631--668.

\bibitem{MR2245536}
\bysame, \emph{On the variety of {L}agrangian subalgebras. {II}}, Ann. Sci.
  \'Ecole Norm. Sup. (4) \textbf{39} (2006), no.~2, 347--379.

\bibitem{MR2372206}
\bysame, \emph{Poisson geometry of the {G}rothendieck resolution of a complex
  semisimple group}, Mosc. Math. J. \textbf{7} (2007), no.~4, 613--642.

\bibitem{MR1726784}
S.~Evens, J.-H. Lu, and A.~Weinstein, \emph{Transverse measures, the modular
  class and a cohomology pairing for {L}ie algebroids}, Quart. J. Math. Oxford
  Ser. (2) \textbf{50} (1999), no.~200, 417--436.

\bibitem{Gualtieri2007}
M.~Gualtieri, \emph{Branes on {P}oisson varieties},  (2007),
  \mbox{\arxiv{0710.2719}}.

\bibitem{Gualtieri2007a}
\bysame, \emph{Generalized complex geometry},  (2007),
  \mbox{\arxiv{math/0401221}}.

\bibitem{MR2013140}
N.~Hitchin, \emph{Generalized {C}alabi-{Y}au manifolds}, Q. J. Math.
  \textbf{54} (2003), no.~3, 281--308.

\bibitem{MR2217300}
\bysame, \emph{Instantons, {P}oisson structures and generalized {K}\"ahler
  geometry}, Comm. Math. Phys. \textbf{265} (2006), no.~1, 131--164.

\bibitem{MR1396600}
J.~Huebschmann, \emph{Poisson geometry of certain moduli spaces}, The
  Proceedings of the Winter School ``Geometry and Physics'' (Srn\'\i, 1994),
  no.~39, 1996, 15--35.

\bibitem{MR1696093}
\bysame, \emph{Duality for {L}ie-{R}inehart algebras and the modular class}, J.
  Reine Angew. Math. \textbf{510} (1999), 103--159.

\bibitem{MR2093043}
D.~Huybrechts, \emph{Complex geometry}, Universitext, Springer-Verlag, Berlin,
  2005, An introduction.

\bibitem{MR1178029}
K.~Iwasaki, \emph{Fuchsian moduli on a {R}iemann surface---its {P}oisson
  structure and {P}oincar\'e-{L}efschetz duality}, Pacific J. Math.
  \textbf{155} (1992), no.~2, 319--340.

\bibitem{MR2055289}
A.~Kapustin, \emph{Topological strings on noncommutative manifolds}, Int. J.
  Geom. Methods Mod. Phys. \textbf{1} (2004), no.~1-2, 49--81.

\bibitem{MR1449222}
D.~Korotkin and H.~Samtleben, \emph{On the quantization of isomonodromic
  deformations on the torus}, Internat. J. Modern Phys. A \textbf{12} (1997),
  no.~11, 2013--2029.

\bibitem{MR1421686}
Y.~Kosmann-Schwarzbach, \emph{The {L}ie bialgebroid of a {P}oisson-{N}ijenhuis
  manifold}, Lett. Math. Phys. \textbf{38} (1996), no.~4, 421--428.

\bibitem{MR1077465}
Y.~Kosmann-Schwarzbach and F.~Magri, \emph{Poisson-{N}ijenhuis structures},
  Ann. Inst. H. Poincar\'e Phys. Th\'eor. \textbf{53} (1990), no.~1, 35--81.

\bibitem{Kosmann-Schwarzbach2007}
\bysame, \emph{On the modular classes of {P}oisson-{N}ijenhuis manifolds},
  (2007), \mbox{\arxiv{math/0611202}}.

\bibitem{MR1472888}
Z.-J. Liu, A.~Weinstein, and P.~Xu, \emph{Manin triples for {L}ie
  bialgebroids}, J. Differential Geom. \textbf{45} (1997), no.~3, 547--574.

\bibitem{private}
J.-H. Lu, private communication.

\bibitem{MR1430434}
\bysame, \emph{Poisson homogeneous spaces and {L}ie algebroids associated to
  {P}oisson actions}, Duke Math. J. \textbf{86} (1997), no.~2, 261--304.

\bibitem{Mackenzie2007}
K.~C.~H. Mackenzie, \emph{Ehresmann doubles and {D}rinfel'd doubles for {L}ie
  algebroids and {L}ie bialgebroids},  (2007), \mbox{\arxiv{math/0611799}}.

\bibitem{MR1262213}
K.~C.~H. Mackenzie and P.~Xu, \emph{Lie bialgebroids and {P}oisson groupoids},
  Duke Math. J. \textbf{73} (1994), no.~2, 415--452.

\bibitem{MR773513}
F.~Magri and C.~Morosi, \emph{On the reduction theory of the {N}ijenhuis
  operators and its applications to {G}el\cprime fand-{D}iki\u\i\ equations},
  Proceedings of the IUTAM-ISIMM symposium on modern developments in analytical
  mechanics, Vol. II (Torino, 1982), vol. 117, 1983, 599--626.

\bibitem{MR900387}
\bysame, \emph{Old and new results on recursion operators: an algebraic
  approach to {KP} equation}, Topics in soliton theory and exactly solvable
  nonlinear equations (Oberwolfach, 1986), World Sci. Publishing, Singapore,
  1987, 78--96.

\bibitem{MR1460632}
T.~Mokri, \emph{Matched pairs of {L}ie algebroids}, Glasgow Math. J.
  \textbf{39} (1997), no.~2, 167--181.

\bibitem{MR1980616}
R.~Moraru, \emph{Integrable systems associated to a {H}opf surface}, Canad. J.
  Math. \textbf{55} (2003), no.~3, 609--635.

\bibitem{MR2228339}
A.~Pichereau, \emph{Poisson (co)homology and isolated singularities}, J.
  Algebra \textbf{299} (2006), no.~2, 747--777.

\bibitem{MR1465521}
A.~Polishchuk, \emph{Algebraic geometry of {P}oisson brackets}, J. Math. Sci.
  (New York) \textbf{84} (1997), no.~5, 1413--1444.

\bibitem{MR842417}
M.~A. Semenov-Tian-Shansky, \emph{Dressing transformations and {P}oisson group
  actions}, Publ. Res. Inst. Math. Sci. \textbf{21} (1985), no.~6, 1237--1260.

\bibitem{MR2276462}
M.~Sti{\'e}non and P.~Xu, \emph{Poisson quasi-{N}ijenhuis manifolds}, Comm.
  Math. Phys. \textbf{270} (2007), no.~3, 709--725.

\bibitem{MR2104604}
L.~Stolovitch, \emph{Sur les structures de {P}oisson singuli\`eres}, Ergodic
  Theory Dynam. Systems \textbf{24} (2004), no.~5, 1833--1863.

\bibitem{MR1390832}
I.~Vaisman, \emph{Complementary {$2$}-forms of {P}oisson structures},
  Compositio Math. \textbf{101} (1996), no.~1, 55--75.

\bibitem{MR2285039}
A.~Weinstein, \emph{The integration problem for complex {L}ie algebroids}, From
  geometry to quantum mechanics, Progr. Math., vol. 252, Birkh\"auser Boston,
  Boston, MA, 2007, 93--109.

\bibitem{MR1675117}
P.~Xu, \emph{Gerstenhaber algebras and {BV}-algebras in {P}oisson geometry},
  Comm. Math. Phys. \textbf{200} (1999), no.~3, 545--560.

\end{thebibliography}
\end{document}